\documentclass[12pt]{article}
\usepackage{amssymb}
\usepackage{amsmath}
\usepackage{latexsym}

\newtheorem{theorem}{Theorem}[section]
\newtheorem{proposition}[theorem]{Proposition}
\newtheorem{corollary}[theorem]{Corollary}
\newtheorem{lemma}[theorem]{Lemma}

\newtheorem{remark}[theorem]{Remark}

\newcommand{\tpi}{2\pi i}
\newcommand{\vac}{\mathbf{1}}
\newcommand{\Tr}{\mathrm{Tr}\, }
\begin{document}

\title{Free Bosonic Vertex Operator Algebras on Genus Two Riemann Surfaces~II}
\author{Geoffrey Mason\thanks{%
Supported by the NSF and NSA} \\
%EndAName
Department of Mathematics, \\
University of California Santa Cruz, \\
CA 95064, U.S.A. \and Michael P. Tuite\\
%EndAName
School of Mathematics, Statistics and Applied Mathematics, \\
National University of Ireland Galway \\
University Road, Galway, Ireland.}
\date{}
\maketitle

\begin{abstract}
We continue our program to define and study $n$-point correlation
functions for a vertex operator algebra $V$ on a 
higher genus compact Riemann surface obtained by sewing surfaces of lower genus.
Here we consider Riemann surfaces of genus $2$ obtained by attaching
a handle to a torus. 
We obtain closed formulas for the genus two partition function for 
 free bosonic theories and lattice vertex operator algebras $V_L$. 
We prove that the partition function is holomorphic in the sewing
parameters on a given suitable domain and describe its modular properties. 
We also compute the genus two Heisenberg vector $n$-point function and show that the
Virasoro vector one point function satisfies a genus two Ward identity.
We compare our results with those obtained in the companion paper,  when a pair of tori are sewn together, and show
that the partition functions are not compatible in the
neighborhood of a two-tori degeneration point. The \emph{normalized}
partition functions of a lattice theory $V_L$ \emph{are} compatible, each being identified with
the genus two theta function of $L$.
\end{abstract}

\newpage

\tableofcontents

\section{Introduction}
\label{Sect_Intro}
In previous work (\cite{MT1}-\cite{MT4}) we developed the general theory of $n$-point functions
for a Vertex Operator Algebra (VOA) on a compact Riemann surface $\mathcal{S}$ obtained by sewing together two surfaces of lower genus, 
and applied this theory to obtain detailed results in the case that $\mathcal{S}$ is obtained by sewing a pair of  complex tori - 
the so-called $\epsilon$-formalism discussed in the companion paper \cite{MT4}. 
In the present paper we consider in detail the situation when $\mathcal{S}$ results from self-sewing a complex torus, i.e., attaching a handle, which we refer to as the $\rho$-formalism. 
We describe the nature of the resulting $n$-point functions, paying particular attention to the $0$-point function, 
i.e., the genus $2$ \emph{partition function}, in the $\rho$-formalism. 
We find the explicit form of the partition function for the Heisenberg 
free bosonic string and for lattice vertex operator algebras, and show that these functions are holomorphic on the parameter domain defined by the sewing. 
We  study the generating function for genus two Heisenberg  $n$-point functions and show that the Virasoro vector 1-point function satisfies a genus two Ward identity.
Many of these results are analogous to those found in the $\epsilon$-formalism discussed in \cite{MT4} but with significant technical differences.
Finally, we compare the results in the two formalisms, and show that the partition functions (and hence all $n$-point functions) are \emph{incompatible}. We introduce \emph{normalized} partition functions,
and in the case of $V_L$ show that they \emph{are} compatible; in both formalisms the normalized
partition function is the genus two Siegel theta function $\theta^{(2)}_L$.

\medskip 
We now discuss the contents of the paper in more detail.
Our approach to genus two correlation functions in both formalisms
is to define them in terms of genus one data coming from a VOA $V$. 
In Section \ref{Sect_Genus2rho} we review the $\rho$-formalism introduced in  \cite{MT2}.
There, we constructed a genus two surface by self-sewing a torus, and obtained explicit expressions for the genus two normalized 2-form of the second kind $\omega^{(2)}$, 
a basis of normalized holomorphic 1-forms $\nu_{1},\nu_{2}$, and the period matrix $\Omega$, in terms of genus one data. 
In particular,     
we constructed a holomorphic map
\begin{eqnarray}
F^{\rho}:\mathcal{D}^{\rho} &\longrightarrow &\mathbb{H}_2
\notag
\\
(\tau ,w,\rho )& \longmapsto & \Omega(\tau ,w,\rho )\label{Frhodef}
\end{eqnarray}
Here, and below, $\mathbb{H}_{g}\ (g \geq 1)$ is the genus $g$ \emph{Siegel upper half-space}, and
$\mathcal{D}^{\rho }\subseteq \mathbb{H}_{1}\times \mathbb{C}^{2}$ is the domain defined
in terms of data $(\tau ,w,\rho )$ needed to self-sew a torus of modulus
$\tau$. Sewing produces a surface $\mathcal{S}= \mathcal{S}(\tau, w, \rho)$ of genus $2$, and the map $F^{\rho}$ assigns to $\mathcal{S}$ its period matrix. See\ Section~\ref{Subsect_rho} for further details. 
In Section~\ref{Subsect_Graphs} we introduce some diagrammatic techniques which
provide a convenient way of describing $\omega^{(2)}, \nu_{1},\nu_{2}$ and $\Omega $
in the $\rho $-formalism. Similar methods are
employed in Sections 5 and 6 to explicitly compute the genus two partition function for free bosonic 
and lattice VOAs. 

\medskip
Section~\ref{Sect_VOA} consists of a brief review of
relevant background material on VOA theory, with particular attention paid
to the Li-Zamolodchikov or LiZ metric. 
In Section~\ref{Sect_Zg},  motivated by ideas in conformal field
theory (\cite{FS}, \cite{So1}, \cite{So2}, \cite{P}),  we introduce $n$-point functions
(at genus one and two) in the
$\rho$-formalism for a general VOA with nondegenerate LiZ metric.
In particular, the genus two \emph{partition function} $Z_{V}^{(2)}: \mathcal{D}^{\rho}\rightarrow \mathbb{C}$ is formally defined as 
\begin{eqnarray}
Z_{V}^{(2)}(\tau ,w,\rho ) &=&\sum_{n\geq 0}\rho ^{n}\sum_{u\in
V_{[n]}}Z_{V}^{(1)}(\bar{u},u,w,\tau ),  \label{pfrho}
\end{eqnarray}%
where $Z_{V}^{(1)}(\bar{u},u,w,\tau )$ is a  genus
one $2$-point function and $\bar{u}$ is the LiZ
metric dual of $u$. 
 In Section 4.1 we consider an example of
self-sewing a sphere (Theorem \ref{Theorem_Z1_q_rho}), while in Section 4.2 we show (Theorem 4.4) that a particular degeneration of the genus $2$ partition function of a VOA $V$ can be described in terms of genus $1$ data. Of particular interest here
is the unexpected r\^{o}le played by the \emph{Catalan series}.

\medskip In Sections~\ref{Sect_HVOA} and \ref{Sect_LVOA} we consider in detail the case of the Heisenberg free bosonic theory $M^l$
corresponding to $l$ free bosons,
and lattice VOAs $V_L$ associated with a positive-definite even lattice $L$. 
Although (\ref{pfrho}) 
is \emph{a priori} a formal power series in
 $\rho ,w$ and $q =e^{2\pi i\tau}$, we will see that
for these two theories it is  a holomorphic function on 
$\mathcal{D}^{\rho}$. We expect that this result holds in much wider
generality. Although our calculations in these two Sections generally parallel those for the $\epsilon$-formalism (\cite{MT4}),
the $\rho$-formalism is not merely a simple translation.
Several issues require additional
attention, so that the $\rho$-formalism is rather more complicated
than its $\epsilon$- counterpart. This circumstance was already evident in 
\cite{MT2}, and stems in part from the fact that $F^{\rho}$
involves a logarithmic term that is absent in the $\epsilon$-formalism. The moment matrices are 
also more unwieldy, as they must be considered as having
entries which themselves are $2\times 2$ block-matrices with entries which
are elliptic-type functions. 

\medskip
We establish 
(Theorem~\ref{Theorem_Z2_boson_rho})  a fundamental formula
describing $Z_{M}^{(2)}(\tau ,w,\rho )$ as a quotient of the genus one partition function
for $M$ by a certain infinite determinant. This determinant already appeared in
\cite{MT2}, and its holomorphy and nonvanishing in $D^{\rho}$ (loc.\ cit.)\ implies the holomorphy
of  $Z_{M}^{(2)}$.  We also obtain
 a product formula for the  infinite determinant (Theorem~\ref{Theorem_Z2_boson_prod_rho}),
and establish the automorphic properties of $Z_{M^{2}}^{(2)}$ with respect to the action of
a group $\Gamma _{1}\cong SL(\mathbb{Z})$ (Theorem~\ref{Theorem_Z2_rho_G1}) 
that naturally
acts on $D^{\rho}$. In particular, we find that $Z^{(2)}_{M^{24}}$ is a form of weight $-12$
with respect to the action of $\Gamma_1$.
These are the analogs in the $\rho$-formalism of  results obtained
in Section 6 of \cite{MT4} for the genus two partition function of $M$ in the 
$\epsilon$-formalism.

\medskip
In Section 5.3 we calculate some 
genus two $n$-point functions for the rank one Heisenberg VOA $M$,
specifically the $n$-point function for the weight $1$ Heisenberg vector 
and the 1-point function for the Virasoro vector $\tilde\omega$. 
We show that, up to an overall factor of the genus two partition function,  the formal differential forms associated with these  
$n$-point functions are described in terms of the global symmetric $2$-form $\omega^{(2)}$ 
(\cite{TUY}) 
and the genus two projective connection (\cite{Gu}) respectively. 
Once again, these results are analogous to results obtained in \cite{MT4} in the $\epsilon$-formalism.
Parallel calculations, including a genus two Ward identity involving the Virasoro $1$-point function (Proposition \ref{Prop: Lattice Virasoro}), are carried out for lattice theories  in Section 6.2.

\medskip
In Section 6.1 we establish (Theorem \ref{Theorem_Z2_L_rho}) a basic formula for the genus two partition function for lattice theories in the $\rho$-formalism. The result is
\begin{equation}\label{pf1}
Z_{V_{L} }^{(2)}(\tau ,w,\rho ) = Z_{M^{l} }^{(2)}(\tau ,w,\rho ) \theta _{L}^{(2)}(\Omega), 
\end{equation}
where $\theta _{L}^{(2)}(\Omega )$ is the genus two Siegel theta function
attached to $L$ (\cite{Fr}) and $\Omega = F^{\rho}(\tau, w, \rho)$; indeed, (\ref{pf1}) is an identity
of formal power series. The holomorphy and automorphic properties of $Z^{(2)}_{V_L, \rho}$ follow from (\ref{pf1}) and those of $Z^{(2)}_{M^l}$ and $\Theta^{(2)}_L$. 

\medskip
Section 7 is devoted to a \emph{comparison} of genus two $n$-point functions, and especially partition
functions, in the $\epsilon$- and $\rho$-formalisms. After all, at genus two we have at our
disposal not one but two partition functions. 
 In view of the strong formal similarities between  $Z_{M^{l},\epsilon }^{(2)}(\tau _{1},\tau _{2},\epsilon
) $ and $Z_{M^{l},\rho }^{(2)}(\tau ,w,\rho )$, for example, it is natural to ask if they are equal in some sense
\footnote{Here we include an additional subscript of either $\epsilon$ or $\rho$ to distinguish between the two formalisms.}. 
One may ask the same question for lattice and other rational theories. 
In the very special case that
$V$ is holomorphic (i.e., it has a \emph{unique} irreducible module),
one knows (e.g., \cite{TUY}) that the genus $2$ conformal block is $1$%
-dimensional, in which case an identification of the two partition functions
might seem inevitable. On the other hand, the partition functions are defined on
quite different domains, so there is no question of them being literally equal. Indeed, we argue
in Section 7 that  $Z_{M^{l},\epsilon }^{(2)}(\tau _{1},\tau _{2},\epsilon)$ 
and $Z_{M^{l},\rho }^{(2)}(\tau ,w,\rho )$ are \emph{incompatible}, i.e.,
there is \emph{no} sensible way in which they can be identified.

  \medskip
  It is useful to introduce \emph{normalized} partition functions, defined as
\begin{equation*}
\widehat{Z}^{(2)}_{V,\rho}(\tau, w, \rho) := \frac{Z_{V,\rho }^{(2)}(\tau ,w,\rho )}{Z_{M^{l},\rho }^{(2)}(\tau ,w,\rho )}, \ \ 
\widehat{Z}^{(2)}_{V,\epsilon}(\tau_1, \tau_2, \epsilon) := \frac{Z_{V,\epsilon}^{(2)}(\tau _{1},\tau _{2},\epsilon )}{Z_{M^{l},\epsilon }^{(2)}(\tau _{1},\tau _{2},\epsilon )}, 
\end{equation*}
associated to a VOA $V$ of central charge $l$. For $M^l$, the normalized
partition functions are equal to $1$. The relation between the normalized partition functions
for lattice theories $V_L\ (\mbox{rk}\,L = l)$ in the two formalisms can be displayed in the diagram 
\begin{equation}\label{commdiag}
\begin{array}{lllll}
\ \hspace{0cm}D^{\epsilon } & \overset{F^{\epsilon }}{\longrightarrow } & \
\ \mathbb{H}_2 & \overset{F^{\rho }}{\longleftarrow } & D^{\rho }  \\ 
& \overset{\widehat{Z}_{V,\epsilon }^{(2)}}{\searrow } & \ \ \ \downarrow \theta
_{L}^{(2)} & \overset{\widehat{Z}_{V,\rho }^{(2)}}{\swarrow } &   \\ 
&  &  &  &  \\ 
&  & \hspace{0.4cm}\mathbb{C} &  & 
\end{array}%
\end{equation}
That this is a \emph{commuting} diagram combines formula (\ref{pf1}) in the $\rho$-formalism,
and Theorem 14 of \cite{MT4} for the analogous result in the $\epsilon$-formalism.
Thus, the \emph{normalized} partition functions for $V_L$
are \emph{independent of the sewing scheme}. They can be identified, via the
sewing maps $F^{\bullet }$, with \emph{a genus two Siegel
modular form of weight $l/2$}, the Siegel theta function. It
is therefore the normalized partition function(s) which can be identified
with an element of the conformal block, and with each other.

\medskip
The fact that (\ref{commdiag}) commutes is, from our current
perspective, something of a miracle; it is achieved by two quite separate and lengthy calculations,
one for each of the formalisms.
It would obviously be useful to have available a result that provides an \emph{a priori} 
guarantee of this fact.
 Section~\ref{Sect_Remarks} contains a brief further discussion of these issues in the light
of related ideas in string theory and algebraic geometry.

\section{Genus Two Riemann Surface from Self-Sewing a Torus}
\label{Sect_Genus2rho}
In this section we review some relevant results of \cite{MT2}
based on a general sewing formalism due to Yamada \cite{Y}. 
In particular, we review the construction of a genus two Riemann 
surface formed by self-sewing a twice-punctured torus. 
We refer to this sewing scheme as the $\rho$-formalism.  
We discuss the explicit form 
of various genus two structures such as the period matrix $\Omega$.
We also review the convergence and holomorphy of an infinite 
determinant that naturally arises later on.
An alternative genus two surface formed by sewing together two tori,
which we refer to as the $\epsilon$-formalism, 
is utilised in the companion paper \cite{MT4}.

\subsection{Some elliptic function theory}
\label{Subsect_Ell}
We begin with the definition of various modular and elliptic functions \cite{MT1}, \cite{MT2}. We define 
\begin{eqnarray}
P_{2}(\tau ,z) &=&\wp (\tau ,z)+E_{2}(\tau )  \notag \\
&=&\frac{1}{z^{2}}+\sum_{k=2}^{\infty }(k-1)E_{k}(\tau )z^{k-2},  \label{P2}
\end{eqnarray}%
where $\tau $ $\in \mathbb{H}_{1}$, the complex upper half-plane and where $%
\wp (\tau ,z)$ is the Weierstrass function 
(with periods $\tpi$ and $\tpi \tau$) and $E_{k}(\tau )=0$ for $k$ odd, and for $k$ even is the Eisenstein series 
\begin{equation}
E_{k}(\tau )=E_{k}(q)=-\frac{B_{k}}{k!}+\frac{2}{(k-1)!}\sum_{n\geq 1}\sigma
_{k-1}(n)q^{n}.  \label{Eisenk}
\end{equation}%
Here and below, we take $q=\exp (2\pi i\tau )$; $\sigma
_{k-1}(n)=\sum_{d\mid n}d^{k-1}$, and $B_{k}$ is a $k$th Bernoulli number
e.g. \cite{Se}. If $k\geq 4$ then $E_{k}(\tau )$ is a holomorphic modular
form of weight $k$ on $SL(2,\mathbb{Z})$ whereas $E_{2}(\tau )$ is a
quasi-modular form \cite{KZ}, \cite{MT2}. We define $P_{0}(\tau ,z)$, up to a
choice of the logarithmic branch, and $P_{1}(\tau ,z)$ by 
\begin{eqnarray}
P_{0}(\tau ,z) &=&-\log (z)+\sum_{k\geq 2}\frac{1}{k}E_{k}(\tau )z^{k},
\label{P0} \\
P_{1}(\tau ,z) &=&\frac{1}{z}-\sum_{k\geq 2}E_{k}(\tau )z^{k-1},  \label{P1}
\end{eqnarray}%
where $P_{2}=-\frac{d}{dz}P_{1}$ and $P_{1}=-\frac{d}{dz}P_{0}.$ $P_{0}$ is
related to the elliptic prime form $K(\tau ,z)$, by \cite{Mu1} 
\begin{equation}
K(\tau ,z)=\exp (-P_{0}(\tau ,z)).  \label{Primeform}
\end{equation}%
Define elliptic functions $P_{k}(\tau ,z)\,$ for $k\geq 3$ 
\begin{equation}
P_{k}(\tau ,z)=\frac{(-1)^{k-1}}{(k-1)!}\frac{d^{k-1}}{dz^{k-1}}P_{1}(\tau
,z).  \label{Pkdef}
\end{equation}

\noindent Define for $k,l\geq 1$ 
\begin{eqnarray}
C(k,l) &=&C(k,l,\tau )=(-1)^{k+1}\frac{(k+l-1)!}{(k-1)!(l-1)!}E_{k+l}(\tau ),
\label{Ckldef} \\
D(k,l,z) &=&D(k,l,\tau ,z)=(-1)^{k+1}\frac{(k+l-1)!}{(k-1)!(l-1)!}%
P_{k+l}(\tau ,z).  \label{Dkldef}
\end{eqnarray}%
$\,$Note that $C(k,l)=C(l,k)$ and $D(k,l,z)=(-1)^{k+l}D(l,k,z)$. We also
define for $k\geq 1$,

\begin{eqnarray}
C(k,0) &=&C(k,0,\tau )=(-1)^{k+1}E_{k}(\tau ),  \label{C(k,0)} \\
D(k,0,z,\tau ) &=&(-1)^{k+1}P_{k}(z,\tau ).  \label{D(k,0)}
\end{eqnarray}%
The Dedekind eta-function is defined by

\begin{equation}
\eta (\tau )=q^{1/24}\prod_{n=1}^{\infty }(1-q^{n}).  \label{etafun}
\end{equation}

\subsection{The $\rho $-formalism for self-sewing a torus}
\label{Subsect_rho}
Consider a compact Riemann surface $\mathcal{S}$ of
genus $2$ with standard homology basis $a_{1},a_{2},b_{1},b_{2}$.
 There exists two holomorphic 1-forms $\nu_{i}$, $i=1,2$ which we may normalize by 
\cite{FK} 
\begin{equation}
\oint_{a_{i}}\nu_{j}=2\pi i\delta_{ij}.  \label{norm}
\end{equation}%
These forms can also be defined via the unique singular bilinear $2$-form $\omega^{(2)}$, known as the 
\emph{normalized differential of the second kind}. 
It is defined by the following properties \cite{FK}, \cite{Y}: 
\begin{equation}
\omega^{(2)}(x,y)=(\frac{1}{(x-y)^{2}}+\text{regular terms})dxdy
\label{omegag}
\end{equation}%
for local coordinates $x,y$, with normalization 
\begin{equation}
\oint_{a_{i}}\omega^{(2)}(x,\cdot )=0\ \ (1 \leq i \leq 2). \label{nugnorm}
\end{equation}% 
Using the Riemann bilinear relations, one finds that 
\begin{equation}
\nu_{i}(x)=\oint_{b_{i}}\omega^{(2)}(x,\cdot ) \ \ (1 \leq i \leq 2).  \label{nui}
\end{equation}%
The period matrix $\Omega$ is defined by 
\begin{equation}
\Omega_{ij}=\frac{1}{2\pi i}\oint_{b_{i}}\nu_{j}\quad \ \ (1 \leq i, j \leq 2). \label{period}
\end{equation}%
It is well-known that $\Omega \in \mathbb{H}_{2}$, the
Siegel upper half plane.

We now review a general method due to Yamada \cite{Y}, and discussed at
length in \cite{MT2}, for calculating $\omega^{(2)}(x,y)$, $\nu_i(x)$ and $\Omega_{ij} $ on the 
Riemann surface formed by sewing
a handle to an oriented torus $\mathcal{S}=\mathbb{C}/\Lambda$ with
lattice $\Lambda=2\pi i(\mathbb{Z}\tau\oplus \mathbb{Z})\ (\tau \in \mathbb{H}_{1})$. 
Consider discs
centered at $z=0$ and $z=w$ with local coordinates $z_{1}=z$ and $z_{2}=z-w$,
and positive radius $r_{a}<\frac{1}{2}D(q)$ with $1\le a\le 2$. Here, we have introduced the minimal lattice distance
\begin{equation}
D(q )=\min_{(m,n)\neq (0,0)}2\pi|m +n\tau |>0.  \label{Dom1om2}
\end{equation}%
Note that $r_{1},r_{2}$ must be sufficiently small to ensure that the discs do not intersect on $\mathcal{S}$. 
Introduce a complex parameter ${\rho }$ where $|{\rho }|\leq
r_{1}r_{2}$ and excise the discs $\{z_{a},\left\vert z_{a}\right\vert \leq |\rho |/r_{\bar{a}}\}$
 to obtain a twice-punctured torus (illustrated in Fig.~1)
\begin{equation*}
\widehat{\mathcal{S}}=\mathcal{S}\backslash \{z_{a},\left\vert
z_{a}\right\vert \leq |\rho |/r_{\bar{a}}\} \ \ (1 \leq a \leq 2).
\end{equation*}%
\begin{center}
\begin{picture}(200,100)

%left surface

\put(50,52){\qbezier(-6,19)(45,35)(91,19)}% upper

\put(50,48){\qbezier(-6,-18)(45,-35)(96,-17)}%left lower

\put(50,50){\qbezier(25,0)(45,17)(60,0)}%upper
\put(50,50){\qbezier(20,2)(45,-17)(65,2)}%lower

% right annulus centered at (140,50)

\put(140,50){\circle{16}}
\put(140,50){\circle{40}}

% z_2=0label
\put(140,50){\vector(-1,-2){0}}%arrow
\put(50,50){\qbezier(90,0)(100,15)(90,30)}%
\put(140,90){\makebox(0,0){$z_2=0$}}

% line and r2 label
\put(140,50){\line(-1,1){14.1}}
\put(127,55){\makebox(0,0){$r_2$}}

% line and eps/r1 label
\put(140,50){\line(1,0){8}}
\put(145,50){\vector(1,4){0}}%arrow
\put(55,20){\qbezier(90,4)(85,17)(90,30)}%
\put(150,15){\makebox(0,0){$|\epsilon|/r_1$}}

%S label
\put(100,30){\makebox(0,0){$\widehat{\mathcal{S}}$}}

%left annulus centred at (45,50) was 185 minus 140

\put(45,50){\circle{16}}
\put(45,50){\circle{40}}

% z_1=0label
\put(45,50){\vector(1,-2){0}}%arrow
\put(0,50){\qbezier(45,0)(35,15)(45,30)}%
\put(45,90){\makebox(0,0){$z_1=0$}}

% line and r1 label
\put(45,50){\line(-1,-1){14.1}}
\put(33,46){\makebox(0,0){$r_1$}}

% line and eps/r2 label
\put(45,50){\line(1,0){8}}
\put(50,50){\vector(1,4){0}}%arrow
\put(0,20){\qbezier(50,4)(45,17)(50,30)}%
\put(50,15){\makebox(0,0){$ |\epsilon|/r_2$}}

\end{picture}

{\small Fig.~1 Self-sewing a torus}
\end{center}
 Here, and below, we use the convention
\begin{equation}
\overline{1}=2,\qquad\overline{2}=1.  \label{bar}
\end{equation}
Define annular regions $\mathcal{A}_{a}=\{z_{a},|{\rho }|r_{%
\bar{a}}^{-1}\leq \left\vert z_{a}\right\vert \leq r_{a}\}\in \widehat{\mathcal{S}}\ (1 \leq a \leq 2)$,
and identify $\mathcal{A}_{1}$ with $\mathcal{A}_{2}$ as a single region via the sewing
relation 
\begin{equation}
z_{1}z_{2}=\rho.  \label{rhosew}
\end{equation}%
The resulting genus two Riemann surface 
(excluding the degeneration point $\rho =0$)
 is parameterized by the domain 
\begin{equation}
\mathcal{D}^{\rho }=\{(\tau ,w,\rho )\in \mathbb{H}_{1}\times \mathbb{%
C\times C}:\  |w-\lambda |>2|\rho |^{1/2}>0\text{ for all }\lambda \in
\Lambda \},  \label{Drho}
\end{equation}%
where the first inequality follows from the requirement that the annuli do not intersect. 
The Riemann surface inherits the genus one homology
basis $a_{1},b_{1}$. The cycle $a_{2}$ is defined to
be the anti-clockwise contour surrounding the puncture at $w$, and 
$b_{2}$ is a path between identified points $z_{1}=z_{0}$ to 
$z_{2}=  \rho/z_{0}$ for some $z_{0}\in\mathcal{A}_{1}$.

\medskip
$\omega^{(2)}, \nu_{i}$ and $(\Omega_{ij})$, are expressed as a functions of 
$(\tau ,w,\rho )\in \mathcal{D}^{\rho }$ in terms of an infinite matrix of $2 \times 2$ blocks
$R(\tau ,w,\rho )=(R(k,l,\tau ,w,\rho )) \ (k,l\geq 1)$ where \cite{MT2} 
\begin{equation}
R(k,l,\tau ,w,\rho )=-\frac{\rho ^{(k+l)/2}}{\sqrt{kl}}
\left(
\begin{array}{cc}
D(k,l,\tau ,w) & C(k,l,\tau ) \\ 
C(k,l,\tau ) & D(l,k,\tau ,w)%
\end{array}%
\right),  \label{Rdef}
\end{equation}%
for $C,D$ of \eqref{Ckldef} and \eqref{Dkldef}.
Note that $R_{ab}(k,l)=R_{\bar{a}\bar{b}}(l,k)\ (1 \leq a, b \leq 2)$. 
$I-R$ and $\det (I-R)$ play a central r\^{o}le in our discussion, where
 $I$ denotes the doubly-indexed identity matrix and $\det (I-R)$ is defined by
 \begin{equation}
\log \det (I-R) =
\Tr\log (I-R)=-\sum_{n\ge 1}\frac{1	}{n}\Tr R^n.
\label{logdetR}
\end{equation} 
In particular (op. cit.  , Proposition 6 and Theorem 7)
\begin{theorem}\label{Theorem_R} We have\\
(a)
\begin{equation}
(I-R)^{-1}=\sum_{n\geq 0}R^{n}  \label{1_minus_R}
\end{equation}%
is convergent in $\mathcal{D}^{\rho }$.\\
(b) $\det (I-R)$ is nonvanishing and holomorphic in $\mathcal{D}^{\rho }$. $\square $
\end{theorem}

We define a set of 1-forms on $\widehat{\mathcal{S}}$  given by
\begin{eqnarray}
a_{1}(k,x) &=&a_{1}(k,x,\tau ,\rho )=\sqrt{k}\rho ^{k/2}P_{k+1}(\tau ,x)dx,
\notag
\\
a_{2}(k,x) &=&a_{2}(k,x,\tau ,\rho )=a_{1}(k,x-w),
\label{aforms}
\end{eqnarray}%
indexed by integers $k\ge 1$. 
We also define the infinite row vector $a(x)=(a_{a}(k,x))$ and
infinite column vector $\overline{a}(x)^{T}=(a_{\bar{a}}(k,x))^{T}\ (k \geq 1, 1 \leq a \leq 2)$.
We find (op. cit., Lemma 11, Proposition 6 and Theorem 9):
\begin{theorem}
\label{theorem_om2} 
\begin{equation}
\omega^{(2)}(x,y)=\omega^{(1)}(x,y)
-a(x)(I-R)^{-1}\overline{a}(y)^{T},
\label{omg2}
\end{equation}
where $\omega^{(1)}(x,y)=P_{2}(\tau, x-y)dxdy$.\quad $\square$

\end{theorem}
Applying (\ref{nui}) then results (op. cit., Lemma 12 and Theorem 9) in 
\begin{theorem}
\label{theorem_nu2} 
\begin{eqnarray}
\nu_{1}(x) 
&=&
dx-\rho^{1/2}\sigma \left((a(x)(I-R)^{-1})(1)\right) \notag
%\\
%&=&
%dx-\rho^{1/2}\sigma \left(((I-R)^{-1}\overline{a}(x)^{T})(1)\right), \label{nu1_rho}
\\
\nu_{2}(x) 
%&=&
%\left(P_{1}(\tau,x-w)-P_{1}(\tau ,x)\right)dx
%-d(I-R)^{-1}\overline{a}(x)^{T} \notag
%\\
&=&
\left(P_{1}(\tau,x-w)-P_{1}(\tau ,x)\right)dx
-a(x)(I-R)^{-1}\overline{d}^{T}. \label{nui_rho}
\end{eqnarray}
$d=(d_{a}(k))$ is a
doubly-indexed infinite row vector\footnote{%
Note that $d$ is denoted by $\beta $ in \cite{MT2}.}%
\begin{eqnarray}
d_{1}(k) 
&=& -\frac{\rho ^{k/2}}{\sqrt{k}}(P_{k}(\tau ,w)-E_{k}(\tau
)), \notag \\
d_{2}(k) 
&=& (-1)^{k}\frac{\rho ^{k/2}}{\sqrt{k}}(P_{k}(\tau ,w)-E_{k}(\tau
)), \label{dk}
\end{eqnarray}%
with $\bar{d}_{a}=d_{\bar{a}}$. $(1)$ refers to the $(k)=(1)$ entry of a row vector and 
$\sigma (M)$ denotes the sum of the entries of a
finite matrix $M$. \quad $\square$
\end{theorem}

$\Omega$ is determined (op. cit., Proposition 11) by \eqref{period} as follows:

\begin{theorem}
\label{Theorem_period_rho} 
There is a holomorphic map%
\begin{eqnarray}
F^{\rho }:\mathcal{D}^{\rho } &\rightarrow &\mathbb{H}_{2},  \notag \\
(\tau ,w,\rho ) &\mapsto &\Omega (\tau ,w,\rho ),  \label{Frhomap}
\end{eqnarray}%
where $\Omega =\Omega (\tau ,w,\rho )$ is given by 
\begin{eqnarray}
2\pi i\Omega _{11} &=&2\pi i\tau -\rho \sigma ((I-R)^{-1}(1,1)),
\label{Om11rho} \\
2\pi i\Omega _{12} &=&w-\rho ^{1/2}\sigma ((d(I-R)^{-1}(1)),  \label{Om12rho}
\\
2\pi i\Omega _{22} &=&\log \left(-\frac{\rho }{K(\tau ,w)^{2}}\right)-d(I-R)^{-1}\bar{d}%
^{T}.  \label{Om22rho}
\end{eqnarray}%
$K$ is the elliptic prime form \eqref{Primeform},  
$(1,1)$ and $(1)$ refer to the $(k,l)=(1,1)$%
, respectively, $(k)=(1)$ entries of an infinite matrix and row vector
respectively. $\sigma (M)$ denotes the sum over the finite indices for a given $2\times 2$ or $1\times 2$ matrix $M$.  
\quad $\square$
\end{theorem}

$\mathcal{D}^{\rho }$ admits an action of the Jacobi group $J=SL(2,%
\mathbb{Z})\ltimes \mathbb{Z}^{2}$ as follows: 
\begin{eqnarray}
(a,b).(\tau ,w,\rho ) &=&(\tau ,w+2\pi ia\tau +2\pi ib,\rho )\quad ((a,b)\in 
\mathbb{Z}^{2}),  \label{ab_rho} \\
\gamma _{1}.(\tau ,w,\rho ) &=&(\frac{a_{1}\tau +b_{1}}{c_{1}\tau +d_{1}},%
\frac{w}{c_{1}\tau +d_{1}},\frac{\rho }{(c_{1}\tau +d_{1})^{2}})\quad
(\gamma _{1}\in \Gamma _{1}),  \label{gam1_rho}
\end{eqnarray}%
with $\Gamma _{1} =
\left\{\left( 
\begin{array}{cc}
a_{1} & b_{1} \\ 
c_{1} & d_{1}%
\end{array}%
\right)\right\} = SL(2,\mathbb{Z})$. 
Due to the branch structure of
the logarithmic term in (\ref{Om22rho}), $F^{\rho }$ is not
equivariant with respect to $J$. 
(Instead one must pass to a simply-connected covering space $\mathcal{\hat{D}}^{\rho }$ 
on which $L=\hat{H}\Gamma _{1}$, a split extension of $SL(2,\mathbb{Z})$ by an integer
Heisenberg group $\hat{H}$, acts. See Section 6.3 of \cite{MT2} for
details). 

\medskip
There is an injection 
$\Gamma _{1}\rightarrow Sp(4,\mathbb{Z})$ defined by 
\begin{equation}
\left( 
\begin{array}{cc}
a_{1} & b_{1} \\ 
c_{1} & d_{1}%
\end{array}%
\right)\mapsto  \left( 
\begin{array}{cccc}
a_{1} & 0 & b_{1} & 0 \\ 
0 & 1 & 0 & 0 \\ 
c_{1} & 0 & d_{1} & 0 \\ 
0 & 0 & 0 & 1%
\end{array}%
\right), \label{Gamma1}
\end{equation}%
through which $\Gamma_1$ acts on
$\mathbb{H}_{2}$ 
by the standard action 
\begin{equation}
\gamma .\Omega {=(A\Omega +B)(C\Omega +D)^{-1}} \ \left(\gamma =\left( 
\begin{array}{ll}
A & B \\ 
C & D%
\end{array}%
\right) \in Sp(4,\mathbb{Z})\right).  \label{eq: modtrans}
\end{equation}%
We then have (op. cit., Theorem 11, Corollary 2)
\begin{theorem}
\label{TheoremG1equiv}  $F^{\rho }$ is
equivariant with respect to the action of $\Gamma _{1}$, i.e.\ there is a
commutative diagram for $\gamma _{1}\in \Gamma _{1}$, 
\begin{equation*}
\begin{array}{ccc}
\mathcal{D}^{\rho } & \overset{F^{\rho }}{\rightarrow } & \mathbb{H}_{2} \\ 
\gamma _{1}\downarrow &  & \downarrow \gamma _{1} \\ 
\mathcal{D}^{\rho } & \overset{F^{\rho }}{\rightarrow } & \mathbb{H}_{2}%
\end{array}
\end{equation*}%
 \quad $\square $
\end{theorem}

\subsection{Graphical expansions}
\label{Subsect_Graphs}
We describe a graphical approach for describing the 
expressions for $\omega,\nu_{i},\Omega _{ij}$ reviewed above. 
These also play an important r\^{o}le in the analysis of genus two partition functions
for the Heisenberg vertex operator algebra. A similar approach is described in \cite{MT4} 
suitable for the $\epsilon$-sewing scheme.
Here we  introduce \emph{doubly-indexed} cycles construed as (clockwise) oriented, labelled polygons $L$
with $n$ nodes for some integer $n\geq 1$, nodes being labelled by a pair of
integers $k,a$ where $k\geq 1$ and $a\in \{1,2\}$. 
Thus, a typical doubly-indexed cycle looks as follows:

\begin{center}
\begin{picture}(250,80)

%\graphpaper(0,0)(250,100)

\put(100,50){\line(1,2){10}}
\put(82,50){\makebox(0,0){$k_1,a_1$}}
\put(100,50){\circle*{4}}
\put(107,63){\vector(1,2){0}}%arrow

\put(110,70){\line(1,0){20}}
\put(100,78){\makebox(0,0){$k_2,a_2$}}
\put(110,70){\circle*{4}}
\put(122,70){\vector(1,0){0}}%arrow

\put(130,70){\line(1,-2){10}}
\put(145,78){\makebox(0,0){$k_3,a_3$}}
\put(130,70){\circle*{4}}
\put(137,57){\vector(1,-2){0}}%arrow

\put(140,50){\line(-1,-2){10}}
\put(160,50){\makebox(0,0){$k_4,a_4$}}
\put(140,50){\circle*{4}}
\put(134,38){\vector(-1,-2){0}}%arrow

\put(110,30){\line(1,0){20}}
\put(145,20){\makebox(0,0){$k_5,a_5$}}
\put(130,30){\circle*{4}}
\put(117,30){\vector(-1,0){0}}%arrow

\put(100,50){\line(1,-2){10}}
\put(100,20){\makebox(0,0){$k_6,a_6$}}
\put(110,30){\circle*{4}}
\put(104,42){\vector(-1,2){0}}%arrow

\end{picture}

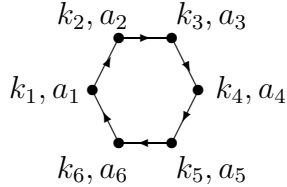
{\small Fig.~2 Doubly-Indexed Cycle}
\end{center}

We define a weight function\footnote{denoted by $\omega$ in Section~6.2 of \cite{MT2}}
 $\zeta $ with values in the ring of
elliptic functions and quasi-modular forms 
$\mathbb{C}[P_{2}(\tau,w),P_{3}(\tau ,w),E_{2}(\tau ),E_{4}(\tau ),E_{6}(\tau )]$ as follows: 
if $L$ is a doubly-indexed cycle then $L$ has edges $E$ labelled as 
$\overset{k,a}{\bullet }\rightarrow \overset{l,b}{\bullet }$, and we set 
\begin{equation}
\zeta (E)=R_{ab}(k,l,\tau ,w,{\rho }),
\label{zetaE}
\end{equation}%
with $R_{ab}(k,l)$ as in (\ref{Rdef}) and 
\begin{equation*}
\zeta (L)=\prod \zeta (E),
\end{equation*}%
where the product is taken over all edges of $L$. 

\medskip We also introduce \emph{doubly-indexed necklaces}
$\mathcal{N}=\{N\} $. These are connected graphs with $n\geq 2$ nodes, $(n-2)$ of which
have valency $2$ and two of which have valency $1$ together with an
orientation, say from left to right, on the edges. 
In this case, each vertex
carries two integer labels $k,a$ with $k\geq 1$ and $a\in \{1,2\}$. 
We define the degenerate necklace $N_{0}$ to be a single node with no edges,
and set $\zeta (N_{0})=1$.

\medskip We define necklaces with distinguished end nodes labelled $k,a;l,b$
as follows: 
\begin{equation*}
\underset{k,a}{\bullet }\longrightarrow \underset{k_{1},a_{1}}{\bullet }%
\ldots   \underset{k_{2},a_{2}}{\bullet }\longrightarrow \underset{l,b}{%
\bullet }\hspace{10mm}\mbox{(type $k,a;l,b$)}
\end{equation*}%
and set 
\begin{equation}
\mathcal{N}(k,a;l,b)=\{\mbox{isomorphism classes of
necklaces of type}\ k,a;l,b\}. 
\label{Nkl}
\end{equation}%
We define 
\begin{eqnarray}
\zeta(1;1) &=&\sum_{a_{1},a_{2}=1,2}\sum_{N\in \mathcal{N}%
(1,a_{1};1,a_{2})}\zeta (N),  \notag \\
\zeta(d;1) &=&
\sum_{a_{1},a_{2}=1,2}
\sum_{k\geq 1}d_{a_{1}}(k)
\sum_{N\in \mathcal{N}(k,a_{1};1,a_{2})}\zeta (N), \notag  \\
\zeta(d;\bar{d}) &=&
\sum_{a_{1},a_{2}=1,2}\sum_{k,l\geq 1}d_{a_{1}}(k)\bar{d}_{a_{2}}(l)
\sum_{N\in \mathcal{N}(k,a_{1};l,a_{2})}\zeta (N).  \label{om_rho_weights}
\end{eqnarray}%
Then we find

\begin{proposition}
\label{Prop_rhoperiod_graph} (\cite{MT2}, Proposition 12) The period matrix is given by 
\begin{eqnarray*}
2\pi i\Omega _{11} &=&2\pi i\tau -\rho \zeta(1;1), \\
2\pi i\Omega _{12} &=&w-\rho ^{1/2}\zeta(d;1), \\
2\pi i\Omega _{22} &=&\log (-\frac{\rho }{K(\tau ,w)^{2}})
-\zeta(d;\bar{d}).\quad \square
\end{eqnarray*}
\end{proposition}

\medskip

We can similarly obtain necklace graphical expansions for the bilinear form 
$\omega^{(2)}(x,y)$ and the holomorphic one forms $\nu_i(x)$. 
We introduce 
further distinguished valence one nodes labelled by $x\in \hat{\mathcal{S}}$, 
the punctured torus. The set of edges $\{E\}$ is augmented by edges with weights defined by: 
\begin{eqnarray}
\zeta(\overset{x}{\bullet}\longrightarrow\overset{y}{\bullet}%
) &=&\omega^{(1)}(x,y),\notag
\\
\zeta(\overset{x}{\bullet}\longrightarrow\overset{k,a}{\bullet}%
) &=&a_{a}(k,x),\notag
\\
\zeta(\overset{k,a}{\bullet}\longrightarrow\overset{y}{\bullet}) &=&-a_{\bar{a}}(k,y),
  \label{eq: zeta1a2}
\end{eqnarray}
with $\omega^{(1)}(x,y)=P_{2}(\tau, x-y)dxdy$ and for 1-forms \eqref{aforms}.

We also consider doubly-indexed necklaces where one or
both end points are $x,y$-labeled nodes. We thus define for $x,y\in \hat{\mathcal{S}}$ 
two isomorphism
classes of oriented doubly-indexed necklaces denoted by  $\mathcal{N}(x;y)$,
and  $\mathcal{N}(x;k,a)$ with the following
respective typical configurations 
\begin{eqnarray}
\{\overset{x}{\bullet }\longrightarrow \overset{k_{1},a_{1}}{\bullet }%
\ldots \overset{k_{2},a_{2}}{\bullet}\longrightarrow\overset{y}{\bullet }%
\},    \label{Nxydef}\\
\{\overset{x}{\bullet }\longrightarrow \overset{k_{1},a_{1}}{\bullet }%
\ldots \overset{k_{2},a_{2}}{\bullet}\longrightarrow\overset{k,a}{\bullet }%
\}.  \label{Nx1def} 
\end{eqnarray}
Furthermore, we define the weights
\begin{eqnarray}
\zeta(x;y) &=& \sum_{N\in \mathcal{N}(x;y)}\zeta (N),  \notag \\
\zeta(x;1) &=&\sum_{a=1,2}
\sum_{N\in \mathcal{N}(x;1,a)}\zeta (N),
\notag\\
\zeta(x;\bar{d}) &=&\sum_{a=1,2}\sum_{k\geq 1}
\sum_{N\in \mathcal{N}(x;k,a)}\zeta (N)\bar{d}_{a}(k).
\label{nu_rho_weights}
\end{eqnarray}%
Comparing to \eqref{omg2} and \eqref{nui_rho} we find the following
graphical expansions for the bilinear form $\omega^{(2)}(x,y)$ and the
holomorphic one forms $\nu_i(x)$

\begin{proposition}
\label{Propomeganugraph} 
For $x, y\in \hat{\mathcal{S}}$
\begin{eqnarray}
\omega^{(2)}(x,y)&=&\zeta(x;y), \label{om_graph} \\
\nu_1(x)&=& dx-\rho ^{1/2}\zeta(x;1), \label{nu1_graph}\\
\nu_2(x)&=& (P_{1}(\tau,x-w)-P_{1}(\tau,x))dx-\zeta(x;\bar{d}). \label{nu2_graph}
\end{eqnarray}
\end{proposition}

\section{Vertex Operator Algebras and the Li-Zamolodchikov metric}
\label{Sect_VOA}
\subsection{Vertex operator algebras}
\label{Subsect_VOA}
We review some relevant aspects of vertex operator algebras (\cite{FHL}, \cite%
{FLM}, \cite{Ka}, \cite{LL}, \cite{MN}, \cite{MT5}). A vertex operator algebra (VOA) is
a quadruple $(V,Y,\mathbf{1},\omega )$ consisting of a $\mathbb{Z}$-graded
complex vector space $V=\bigoplus_{n\in \mathbb{Z}}V_{n}$, a linear map $%
Y:V\rightarrow (\mathrm{End}\,V)[[z,z^{-1}]]$, for formal parameter $z$, and a
pair of distinguished vectors (states), the vacuum $\vac \in V_{0}$ ,
and the conformal vector $\omega \in V_{2}$. For each state $v\in V$ the
image under the $Y$ map is the vertex operator

\begin{equation}
Y(v,z)=\sum_{n\in \mathbb{Z}}v(n)z^{-n-1},  \label{Ydefn}
\end{equation}%
with modes $v(n)\in \mathrm{End}\,V$ where $\mathrm{Res}_{z=0}z^{-1}Y(v,z)%
\vac =v(-1)\vac =v$. Vertex operators satisfy the Jacobi identity
or equivalently, operator locality or Borcherds's identity for the modes
(loc. cit.).

\medskip
The vertex operator for the conformal vector $\omega $ is defined as
\begin{equation*}
Y(w,z)=\sum_{n\in \mathbb{Z}}L(n)z^{-n-2}.
\end{equation*}%
The modes $L(n)$ satisfy the Virasoro algebra of central charge $c$:
\begin{equation*}
\lbrack L(m),L(n)]=(m-n)L(m+n)+(m^{3}-m)\frac{c}{12}\delta _{m,-n}.
\end{equation*}
We define the homogeneous space of weight $k$ to be $V_{k}=\{v\in
V|L(0)v=kv\}$ where we write $wt(v)=k$ for $v$ in $V_{k}$ . Then as an
operator on $V$ we have
\begin{equation*}
v(n):V_{m}\rightarrow V_{m+k-n-1}.
\end{equation*}%
In particular, the \textit{zero mode} $o(v)=v(wt(v)-1)$ is a linear operator
on $V_{m}$. A state $v$ is said to be \textit{quasi-primary} if $L(1)v=0$
and \textit{primary} if additionally $L(2)v=0$.

\medskip
The subalgebra $\{L(-1),L(0),L(1)\}$ generates a natural action on vertex
operators associated with $SL(2,\mathbb{C})$ M\"{o}bius transformations 
(\cite{B}, \cite{DGM}, \cite{FHL}, \cite{Ka}). In particular, we note
the inversion $z\mapsto 1/z$, for which 
\begin{equation}
Y(v,z)\mapsto Y^{\dagger }(v,z)=Y(e^{zL(1)}(-\frac{1}{z^{2}})^{L(0)}v,\frac{1%
}{z}).  \label{eq: adj op}
\end{equation}%
$Y^{\dagger }(v,z)$ is the \emph{adjoint} vertex operator \cite{FHL}. 
%Under
%the dilatation $z\mapsto az$ we have 
%\begin{equation}
%Y(v,z)\mapsto a^{L(0)}Y(v,z)a^{-L(0)}=Y(a^{L(0)}v,az).  \label{Y_D}
%\end{equation}%

\medskip
 We consider in particular  the Heisenberg free boson VOA and
lattice VOAs. Consider an $l$-dimensional complex vector space (i.e.,
abelian Lie algebra) $\mathfrak{H}$ equipped with a non-degenerate,
symmetric, bilinear form $(\ ,)$ and a distinguished orthonormal basis $%
a_{1},a_{2},\ldots  a_{l}$. The corresponding affine Lie algebra is the
Heisenberg Lie algebra $\mathfrak{\hat{H}}=\mathfrak{H}\otimes \mathbb{C}%
[t,t^{-1}]\oplus \mathbb{C}k$ with brackets $[k,\mathfrak{\hat{H}}]=0$ and
\begin{equation}
\lbrack a_{i}\otimes t^{m},a_{j}\otimes t^{n}]=m\delta _{i,j}\delta _{m,-n}k.
\label{Fockbracket}
\end{equation}%
Corresponding to an element $\lambda $ in the dual space $\mathfrak{H}^{\ast
}$ we consider the Fock space defined by the induced (Verma) module 
\begin{equation*}
M^{(\lambda )}=U(\mathfrak{\hat{H}})\otimes _{U(\mathfrak{H}\otimes \mathbb{C%
}[t]\oplus \mathbb{C}k)}\mathbb{C},
\end{equation*}%
where $\mathbb{C}$ is the $1$-dimensional space annihilated by $\mathfrak{H}%
\otimes t\mathbb{C}[t]$ and on which $k$ acts as the identity and $\mathfrak{%
H}\otimes t^{0}$ via the character $\lambda $; $U$ denotes the universal
enveloping algebra. There is a canonical identification of linear spaces
\begin{equation*}
M^{(\lambda )}=S(\mathfrak{H}\otimes t^{-1}\mathbb{C}[t^{-1}]),
\end{equation*}%
where $S$ denotes the (graded) symmetric algebra. The Heisenberg free boson
VOA $M^{l}$ corresponds to the case $\lambda =0$. The Fock states
\begin{equation}
v=a_{1}(-1)^{e_{1}}.a_{1}(-2)^{e_{2}}\ldots  .a_{1}(-n)^{e_{n}}\ldots  .a_{l}(-1)^{f_{1}}.a_{l}(-2)^{f_{2}}\ldots  a_{l}(-p)^{f_{p}}.%
\mathbf{1,}  \label{Fockstate}
\end{equation}%
for non-negative integers $e_{i},\ldots  ,f_{j}$ form a\textit{\ }basis of $M^{l}$%
. The vacuum $\vac $ is canonically identified with the identity of $%
M_{0}=\mathbb{C}$, while the weight 1 subspace $M_{1}$ may be naturally
identified with $\mathfrak{H}$. $M^{l}$ is a simple VOA of central charge $l$.

\medskip Next we consider the case of a lattice vertex operator algebra $%
V_{L}$ associated to a positive-definite even lattice $L$ (cf. \cite{B}, 
\cite{FLM}). Thus $L$ is a free abelian group of rank $l$ equipped with a
positive definite, integral bilinear form $(\ ,):L\otimes L\rightarrow 
\mathbb{Z}$ such that $(\alpha ,\alpha )$ is even for $\alpha \in L$. Let $%
\mathfrak{H}$ be the space $\mathbb{C}\otimes _{\mathbb{Z}}L$ equipped with
the $\mathbb{C}$-linear extension of $(\ ,)$ to $\mathfrak{H}\otimes 
\mathfrak{H}$ and let $M^{l}$ be the corresponding Heisenberg VOA. The Fock
space of the lattice theory may be described by the linear space 
\begin{equation}
V_{L}=M^{l}\otimes \mathbb{C}[L]=\sum_{\alpha \in L}M^{l}\otimes e^{\alpha },
\label{VLdefn}
\end{equation}%
where $\mathbb{C}[L]$ denotes the group algebra of $L$ with canonical basis $%
e^{\alpha }$, $\alpha \in L$. $M^{l}$ may be identified with the subspace $%
M^{l}\otimes e^{0}$ of $V_{L}$, in which case $M^{l}$ is a subVOA of $V_{L}$
and the rightmost equation of (\ref{VLdefn}) then displays the decomposition
of $V_{L}$ into irreducible $M^{l}$-modules. $V_{L}$ is a simple VOA of
central charge $l$. Each $\vac \otimes e^{\alpha }\in V_{L}$ is a
primary state of weight $\frac{1}{2}(\alpha ,\alpha )$ with vertex operator
(loc. cit.) 
\begin{eqnarray}
Y(\vac \otimes e^{\alpha },z) &=&
Y_{-}(\alpha,z)Y_{+}(\alpha ,z)e^{\alpha }z^{\alpha },  \notag \\
Y_{\pm }(\alpha ,z) &=&\exp (\mp \sum_{n>0}\frac{%
\alpha (\pm n)}{n}z^{\mp n}).  \label{Yealpha}
\end{eqnarray}%
The operators $e^{\alpha }\in \mathbb{C}[L]$ obey 
\begin{equation}
e^{\alpha }e^{\beta }=\epsilon (\alpha ,\beta )e^{\alpha +\beta }
\label{eps_cocycle}
\end{equation}%
for a bilinear $2$-cocycle $\epsilon (\alpha ,\beta )$ satisfying $\epsilon (\alpha
,\beta )\epsilon (\beta ,\alpha )=(-1)^{(\alpha ,\beta )}$.  

\subsection{The Li-Zamolodchikov metric}
\label{Subsect_LiZ}
A bilinear form $\langle \ ,\rangle :V\times V{\longrightarrow }\,\mathbb{C}$
is called \emph{invariant} in case the following identity holds for all $%
a,b,c\in V$ (\cite{FHL}): 
\begin{equation}
\langle Y(a,z)b,c\rangle =\langle b,Y^{\dagger }(a,z)c\rangle,
\label{eq: inv bil form}
\end{equation}%
with $Y^{\dagger }(a,z)$ the adjoint operator (\ref{eq: adj op}).
If $V_{0}=\mathbb{C}\vac $ and $V$ is self-dual (i.e. $V$ is isomorphic to the
contragredient module $V^{\prime }$ as a $V$-module) then $V$ has a unique
non-zero invariant bilinear form up to scalar \cite{Li}. 
Note that $\langle \ ,\rangle 
$ is necessarily symmetric by a theorem of \cite{FHL}. Furthermore, if $V$ is
simple then such a form is necessarily non-degenerate. All of the VOAs that
occur in this paper satisfy these conditions, so that normalizing $\langle \vac ,\vac \rangle =1$ implies that $\langle \ ,\rangle $ is
unique. We refer to such a bilinear form as the 
\emph{Li-Zamolodchikov metric} on $V$, or LiZ-metric for short \cite{MT4}.
We also note that the LiZ-metric is multiplicative over tensor products in the sense that LiZ metric of the tensor
product $V_{1}\otimes V_{2}$ of a pair of simple VOAs satisfying the
above conditions is by uniqueness, the tensor product of the LiZ metrics on $V_{1}$ and $V_{2}$.

\medskip
For a quasi-primary vector $a$ of weight $wt(a)$, the component form of \eqref{eq: inv bil form} becomes 
\begin{equation}
\langle a(n)b,c\rangle =(-1)^{wt(a)}\langle b,a(2wt(a)-n-2)c\rangle .
\label{eq: qp bil form}
\end{equation}%
In particular, for the conformal vector $\omega $ we obtain 
\begin{equation}
\langle L(n)b,c\rangle =\langle b,L(-n)c\rangle .
\label{eq: conf vec bil form}
\end{equation}%
Taking $n=0$, it follows that the homogeneous
spaces $V_{n}$ and $V_{m}$ are orthogonal if $n\not=m$. 

\medskip Consider the rank one Heisenberg VOA $M = M^1$
generated by a weight one state $a$ with $(a,a)=1$. Then $\langle a,a\rangle
=-\langle \vac ,a(1)a(-1)\vac \rangle =-1$. Using (\ref{Fockbracket}),
 it is straightforward to verify that the Fock
basis \eqref{Fockstate} is orthogonal with respect to the
LiZ-metric and 
\begin{equation}
\langle v,v\rangle =\prod_{1\leq i\leq p}(-i)^{e_{i}}e_{i}!.
\label{eq: inner prod}
\end{equation}%
This result generalizes in an obvious way to the rank $l$ free boson VOA
$M^{l}$ because the LiZ metric is multiplicative over tensor products.

\medskip We consider next the lattice vertex operator algebra $V_{L}$ for a
positive-definite even lattice $L$. We take as our Fock basis the states $%
\{v\otimes e^{\alpha }\}$ where $v$ is as in (\ref{Fockstate}) and $\alpha $
ranges over the elements of $L$.

\begin{lemma}
\label{Lemma_LiZ_lattice} If $u,v\in M^{l}$ and $\alpha ,\beta \in L$, then 
\begin{eqnarray*}
\langle u\otimes e^{\alpha },v\otimes e^{\beta }\rangle &=&\langle
u,v\rangle \langle \vac \otimes e^{\alpha },\vac \otimes e^{\beta
}\rangle \\
&=&(-1)^{\frac{1}{2}(\alpha ,\alpha )}\epsilon (\alpha ,-\alpha )\langle
u,v\rangle \delta _{\alpha ,-\beta }.
\end{eqnarray*}
\end{lemma}
\noindent \textbf{Proof.} It follows by successive applications of (\ref{eq:
qp bil form}) that the first equality in the Lemma is true, and that it is
therefore enough to prove it in the case that $u=v=\vac $. We identify
the primary vector $\vac \otimes e^{\alpha }$ with $e^{\alpha }$ in the
following. Then $\langle
e^{\alpha },e^{\beta }\rangle =\langle e^{\alpha}(-1) \vac ,e^{\beta }\rangle$ is given by 
\begin{eqnarray*}
&&(-1)^{\frac{1}{2}(\alpha ,\alpha )}
\langle 
\vac,
e^{\alpha}((\alpha ,\alpha )-1)e^{\beta }
\rangle \\
&&= (-1)^{\frac{1}{2}(\alpha ,\alpha )}\,
\mathrm{Res}_{z=0}z^{(\alpha ,\alpha )-1}
\langle \vac,
Y(e^{\alpha },z)e^{\beta }
\rangle \\
&&=(-1)^{\frac{1}{2}(\alpha ,\alpha )}
\epsilon (\alpha ,\beta )\,
\mathrm{Res}_{z=0}z^{(\alpha ,\beta )+(\alpha ,\alpha )-1}
\langle 
\vac,
Y_{-}(\alpha ,z).e^{\alpha +\beta }
\rangle.
\end{eqnarray*}
Unless $\alpha +\beta =0$, all states to the left inside the bracket $%
\langle \ ,\rangle $ on the previous line have positive weight, hence are
orthogonal to $\vac $. So $\langle e^{\alpha },e^{\beta }\rangle =0$ if 
$\alpha +\beta \not=0$. In the contrary case, the exponential operator
acting on the vacuum yields just the vacuum itself among weight zero states,
and we get $\langle e^{\alpha },e^{-\alpha }\rangle =(-1)^{\frac{1}{2}%
(\alpha ,\alpha )}\epsilon (\alpha ,-\alpha )$ in this case. \quad $\square $

\begin{corollary}
\label{Corollary_cocycle_choice}We may choose the cocycle so that $\epsilon
(\alpha ,-\alpha )=(-1)^{\frac{1}{2}(\alpha ,\alpha )}$ (cf.\ (\ref{A.1}) in the Appendix). In
this case, we have 
\begin{equation}
\langle u\otimes e^{\alpha },v\otimes e^{\beta }\rangle =\langle u,v\rangle
\delta _{\alpha ,-\beta }.  \label{LiZ_lattice}
\end{equation}
\end{corollary}

\section{Partition and $n$-Point Functions for Vertex Operator Algebras on a Genus Two Riemann Surface}\label{Sect_Zg}

In this section we consider the partition and $n$-point functions for a VOA
on Riemann surface of genus one or two, formed by attaching a handle to a surface of lower genus.
We assume that $V$ has a non-degenerate LiZ metric $\langle \ ,\rangle $. Then for any $V$
basis $\{u^{(a)}\}$, we may define the \emph{dual basis }$\{\bar{u}^{(a)}\}$
with respect to the LiZ metric where 
\begin{equation}
\langle u^{(a)},\bar{u}^{(b)}\rangle =\delta_{ab}.  \label{LiZdual}
\end{equation}

\subsection{Genus one}
 It is instructive to first consider an alternative approach to defining the genus one partition function.  
In order to define $n$-point correlation functions on a torus, 
Zhu introduced (\cite{Z}) a second VOA $(V,Y[\,,\,],\vac ,\tilde{\omega})$ isomorphic to $(V,Y(\,,\,),\vac ,\omega )$ with
vertex operators 
\begin{equation}
Y[v,z]=\sum_{n\in \mathbb{Z}}v[n]z^{-n-1}=Y(q_{z}^{L(0)}v,q_{z}-1),
\label{Ysquare}
\end{equation}%
and conformal vector $\tilde{\omega}=$ $\omega -\frac{c}{24}\vac $. Let 
\begin{equation}
Y[\tilde{\omega},z]=\sum_{n\in \mathbb{Z}}L[n]z^{-n-2} ,  \label{Ywtilde}
\end{equation}%
and write $wt[v]=k$ if $L[0]v=kv$, $V_{[k]}=\{v\in V|wt[v]=k\}$. Similarly, we define the square bracket LiZ metric 
$\langle \ ,\rangle _{\mathrm{sq}}$ which is invariant with respect to the
square bracket adjoint.

\medskip The (genus one) $1$-point function is now defined as
\begin{equation}
Z_{V}^{(1)}(v,\tau )=\Tr_{V}\left(o(v)q^{L(0)-c/24}\right).
\label{Z1_1_pt}
\end{equation}%
An $n$-point function can be expressed in terms
of $1$-point functions (\cite{MT1}, Lemma 3.1) as follows: 
\begin{eqnarray}
&&Z_{V}^{(1)}(v_{1},z_{1};\ldots   v_{n},z_{n};\tau )  \notag \\
&=&Z_{V}^{(1)}(Y[v_{1},z_{1}]\ldots   Y[v_{n-1},z_{n-1}]Y[v_{n},z_{n}]\mathbf{1%
},\tau )  \label{Z1Ysq1} \\
&=&Z_{V}^{(1)}(Y[v_{1},z_{1n}]\ldots   Y[v_{n-1},z_{n-1n}]v_{n},\tau ),
\label{Z1Ysq2}
\end{eqnarray}%
where $z_{in}=z_{i}-z_{n}\ (1 \leq i \leq n-1)$. In particular, $Z_{V}^{(1)}(v_{1},z_{1};v_{2},z_{2};\tau )$
depends only one $z_{12}$, and we denote this $2$-point function by 
\begin{eqnarray}
Z_{V}^{(1)}(v_{1},v_{2},z_{12},\tau )
&=&Z_{V}^{(1)}(v_{1},z_{1};v_{2},z_{2};\tau )  \notag \\
&=&\Tr_{V}\left(o(Y[v_{1},z_{12}]v_{2})q^{L(0)}\right).  \label{Z1_2pt}
\end{eqnarray}

Now consider a torus obtained by self-sewing a
Riemann sphere with punctures located at the origin and an arbitrary point $w$ on the complex plane (cf.\ \cite{MT2}, Section~5.2.2.).
Choose local coordinates $z_{1}$ in the neighborhood of the origin and $%
z_{2}=z-w$ for $z$ in the neighborhood of $w$. For a complex sewing
parameter $\rho$, identify the annuli $|{\rho }|r_{\bar{a}}^{-1}\leq
\left\vert z_{a}\right\vert \leq r_{a}$ for $1\le a \le 2$ and $|{\rho }|\leq
r_{1}r_{2}$ via the sewing relation 
\begin{equation}
z_{1}z_{2}=\rho .  \label{spheresew2}
\end{equation}%
Define
\begin{eqnarray}\label{morechidef}
 \chi =- \frac{\rho}{w^2}.
 \end{eqnarray}
 Then the annuli do not intersect provided $|\chi |<\frac{1}{4}$, and the torus modular parameter is
\begin{equation}
q=f(\chi ),  \label{qCat}
\end{equation}%
where  $f(\chi )$ the Catalan series 
\begin{equation}
f(\chi )=\frac{1-\sqrt{1-4\chi }}{2\chi }-1=\sum_{n\geq 1}\frac{1}{n}\binom{2n}{n+1}\chi ^{n}.  \label{Catalan}
\end{equation}%
$f=f(\chi )$ satisfies $f=\chi (1+f)^{2}$ and the following identity, which can be proved by induction on $m$:
\begin{lemma}
\label{Lemma_catm}$f(\chi )$ satisfies 
\begin{equation*}
f(\chi )^{m}=\sum_{n\geq m}\frac{m}{n}\binom{2n}{n+m}\chi ^{n}\ \quad (m \geq 1).\quad \square
\end{equation*}
\end{lemma}

\medskip 
We now define the genus one partition function in the sewing scheme
\eqref{spheresew2} by 
\begin{eqnarray}
Z_{V,\rho }^{(1)}(\rho ,w)=  %\notag \\
\sum_{n\geq 0}\rho ^{n}\sum_{u\in V_{n}}\mathrm{Res}_{z_{2}=0}z_{2}^{-1}%
\mathrm{Res}_{z_{1}=0}z_{1}^{-1}\langle \vac ,Y(u,w+z_{2})Y(\bar{u}%
,z_{1})\vac \rangle.  \label{Z1def_rho}
\end{eqnarray}%
This partition function is directly related to $Z_{V}^{(1)}(q)=\Tr_{V}(q^{L(0)-c/24})$
as follows:
\begin{theorem}
\label{Theorem_Z1_q_rho} In the sewing scheme (\ref{spheresew2}), we have 
\begin{equation}
Z_{V,\rho }^{(1)}(\rho ,w)=q^{c/24}Z_{V}^{(1)}(q),  \label{Z1_q_chi}
\end{equation}%
where $q=f(\chi )$ is given by (\ref{qCat}).
\end{theorem}

\noindent \textbf{Proof.}
The summand in (\ref{Z1def_rho}) is 
\begin{eqnarray*}
\langle \mathbf{1,}Y(u,w)\bar{u}\rangle &=&\langle Y^{\dagger }(u,w)\mathbf{1%
},\bar{u}\rangle \\
&=&(-w^{-2})^{n}\langle Y(e^{wL(1)}u,w^{-1})\vac ,\bar{u}\rangle \\
&=&(-w^{-2})^{n}\langle e^{w^{-1}L(-1)}e^{wL(1)}u,\bar{u}\rangle ,
\end{eqnarray*}%
where we have used (\ref{eq: adj op}) and also $Y(v,z)\vac =\exp
(zL(-1))v.$ (See \cite{Ka}\ or \cite{MN} for the latter equality.) Hence we
find that 
\begin{eqnarray*}
Z_{V,\rho }^{(1)}(\rho ,w) &=&\sum_{n\geq 0}(-\frac{\rho }{w^{2}}%
)^{n}\sum_{u\in V_{n}}\langle e^{w^{-1}L(-1)}e^{wL(1)}u,\bar{u}\rangle \\
&=&\sum_{n\geq 0}\chi ^{n}\Tr_{V_{n}}(e^{w^{-1}L(-1)}e^{wL(1)}).
\end{eqnarray*}%
Expanding the exponentials yields 
\begin{equation}
Z_{V,\rho }^{(1)}(\rho ,w)=\Tr_{V}(\sum_{r\geq 0}\frac{L(-1)^{r}L(1)^{r}}{%
(r!)^{2}}\chi ^{L(0)}),  \label{Z1rhotrace}
\end{equation}%
an expression which depends only on $\chi $.

\medskip In order to compute (\ref{Z1rhotrace}) we consider the
quasi-primary decomposition of $V$. Let $Q_{m}=\{v\in V_{m}|L(1)v=0\}$
denote the space of quasiprimary states of weight $m\geq 1$. Then $\dim
Q_{m}=p_{m}-p_{m-1}$ with $p_{m}=\dim V_{m}$. Consider the decomposition of $%
V$ into $L(-1)$-descendents of quasi-primaries 
\begin{equation}
V_{n}=\bigoplus_{m=1}^{n}L(-1)^{n-m}Q_{m}.  \label{V_qprim}
\end{equation}
\begin{lemma}
\label{Lemma_Quasip}Let $v\in Q_{m}$ for $m\geq 1$. For an integer $n\geq m,$
\begin{equation*}
\sum_{r\geq 0}\frac{L(-1)^{r}L(1)^{r}}{(r!)^{2}}L(-1)^{n-m}v=\binom{2n-1}{n-m%
}L(-1)^{n-m}v.
\end{equation*}
\end{lemma}
\noindent \textbf{Proof.} First use induction on $t\geq 0$ to show that 
\begin{eqnarray*}
L(1)L(-1)^{t}v=t(2m+t-1)L(-1)^{t-1}v.
\end{eqnarray*}
Then by induction in $r$ it follows that 
\begin{equation*}
\frac{L(-1)^{r}L(1)^{r}}{(r!)^{2}}L(-1)^{n-m}v=\binom{n-m}{r}\binom{n+m-1}{r}%
L(-1)^{n-m}v.
\end{equation*}%
Hence 
\begin{eqnarray*}
\sum_{r\geq 0}\frac{L(-1)^{r}L(1)^{r}}{(r!)^{2}}L(-1)^{n-m}v &=&\sum_{r\geq
0}^{n-m}\binom{n-m}{r}\binom{n+m-1}{r}L(-1)^{n-m}v, \\
&=&\binom{2n-1}{n-m}L(-1)^{n-m}v,
\end{eqnarray*}%
where the last equality follows from a comparison of the coefficient
of $x^{n-m}$ in the identity $(1+x)^{n-m}(1+x)^{n+m-1}=(1+x)^{2n-1} $.   $\hfill \square $

\medskip Lemma \ref{Lemma_Quasip} \ and (\ref{V_qprim}) imply that for $%
n\geq 1,$ 
\begin{eqnarray*}
\Tr_{V_{n}}(\sum_{r\geq 0}\frac{L(-1)^{r}L(1)^{r}}{(r!)^{2}})
&=&\sum_{m=1}^{n}\Tr_{Q_{m}}(\sum_{r\geq 0}\frac{L(-1)^{r}L(1)^{r}}{(r!)^{2}}%
L(-1)^{n-m}) \\
&=&\sum_{m=1}^{n}(p_{m}-p_{m-1})\binom{2n-1}{n-m}.
\end{eqnarray*}%
The coefficient of $p_{m}$ is 
\begin{equation*}
\binom{2n-1}{n-m}-\binom{2n-1}{n-m-1}=\frac{m}{n}\binom{2n}{m+n},
\end{equation*}%
and hence 
\begin{equation*}
\Tr_{V_{n}}(\sum_{r\geq 0}\frac{L(-1)^{r}L(1)^{r}}{(r!)^{2}})=\sum_{m=1}^{n}%
\frac{m}{n}\binom{2n}{m+n}p_{m}.
\end{equation*}%
Using Lemma \ref{Lemma_catm}, we find that 
\begin{eqnarray*}
Z_{V,\rho }^{(1)}(\rho ,w) &=&1+\sum_{n\geq 1}\chi ^{n}\sum_{m=1}^{n}\frac{m%
}{n}\binom{2n}{m+n}p_{m}, \\
&=&1+\sum_{m\geq 1}p_{m}\sum_{n\geq m}\frac{m}{n}\binom{2n}{m+n}\chi ^{n} \\
&=&1+\sum_{m\geq 1}p_{m}(f(\chi ))^{m}
\\&=& \Tr_{V}(f(\chi )^{L(0)}),
\end{eqnarray*}%
and Theorem \ref{Theorem_Z1_q_rho} follows. $\hfill  \square $

\subsection{Genus two}
\label{Subsect_Z2}
We now turn to the case of genus two. Following Section~\ref{Subsect_rho}, we employ the $\rho $-sewing scheme to self-sew a torus $\mathcal{S}$ with
modular parameter $\tau $ via the sewing relation (\ref{rhosew}). 
For $x_{1},\ldots   ,x_{n}\in \mathcal{S}$ with $\left\vert x_{i}\right\vert \geq
|\epsilon |/r_{2}$ and $\left\vert x_{i}-w\right\vert \geq |\epsilon |/r_{1}$%
, we define the genus two $n$-point function in the $\rho $-formalism by 
\begin{gather}
Z_{V }^{(2)}(v_{1},x_{1};\ldots   v_{n},x_{n};\tau ,w,\rho )=  \notag \\
\sum_{r\geq 0}\rho ^{r}\sum_{u\in V_{[r]}}\mathrm{Res}_{z_{1}=0}z_{1}^{-1}%
\mathrm{Res}_{z_{2}=0}z_{2}^{-1}Z_{V}^{(1)}(\bar{u},w+z_{2};v_{1},x_{1};%
\ldots   v_{n},x_{n};u,z_{1};\tau ).  \label{Z2n_pt_rho}
\end{gather}
 In particular, with the notation (\ref{Z1_2pt}), the genus two partition
function is 
\begin{equation}
Z_{V}^{(2)}(\tau ,w,\rho )=\sum_{n\geq 0}\rho ^{n}\sum_{u\in
V_{[n]}}Z_{V}^{(1)}(\bar{u},u,w,\tau ).  \label{Z2_def_rho}
\end{equation}

\medskip Next we consider $Z_{V}^{(2)}(\tau ,w,\rho )$ in the two-tori
degeneration limit. Define, much as in \eqref{morechidef},
\begin{equation}
\chi=-\frac{\rho}{w^2},
\label{eq:chi2}
\end{equation}
where $w$ denotes a point on the torus and $\rho$ is the genus two sewing parameter. Then one finds that two tori degeneration limit is given by $\rho
,w\rightarrow 0$ for fixed $\chi $, where 
\begin{equation}
\Omega \rightarrow \left( 
\begin{array}{cc}
\tau & 0 \\ 
0 & \frac{1}{2\pi i}\log (f(\chi ))
\end{array} \right)  \label{Omrho_degen}
\end{equation}%
and $f(\chi )$ is the Catalan series \eqref{Catalan} (cf.\ \cite{MT2}, Section~6.4). 
\begin{theorem}
\label{Theorem_Degen_pt}
For fixed $|\chi |<\frac{1}{4}$, we have 
\begin{equation*}
\lim_{w,\rho \rightarrow 0}Z_{V }^{(2)}(\tau ,w,\rho
)=f(\chi )^{c/24}Z_{V}^{(1)}(q)Z_{V}^{(1)}(f(\chi )).
\end{equation*}%
\end{theorem}
\noindent \textbf{Proof.} By \eqref{Z1_2pt} we have
\begin{equation*}
Z_{V}^{(1)}(\bar{u},u,w,\tau )=\Tr_{V}(o(Y[\bar{u},w]u)q^{L(0)}).
\end{equation*}%
Using the non-degeneracy of the LiZ metric 
$\langle\ ,\ \rangle_{\mathrm{sq}}$ in the square bracket formalism we obtain 
\begin{equation*}
Y[\bar{u},w]u=\sum_{m\geq 0}\sum_{v\in V_{[m]}}\langle \bar{v},Y[\bar{u}%
,w]u\rangle _{\mathrm{sq}}v.
\end{equation*}%
Arguing much as in the first part of the proof of Theorem \ref{Theorem_Z1_q_rho}, we also find 
\begin{eqnarray*}
\langle \bar{v},Y[\bar{u},w]u\rangle _{\mathrm{sq}} &=&(-w^{-2})^{n}\langle
Y[e^{wL[1]}\bar{u},w^{-1}]\bar{v},u\rangle _{\mathrm{sq}} \\
&=&(-w^{-2})^{n}\langle e^{w^{-1}L[-1]}Y[\bar{v},-w^{-1}]e^{wL[1]}\bar{u}%
,u\rangle _{\mathrm{sq}} \\
&=&(-w^{-2})^{n}\langle E[\bar{v},w]\bar{u},u\rangle _{\mathrm{sq}},
\end{eqnarray*}%
where 
\begin{equation*}
E[\bar{v},w]=\exp (w^{-1}L[-1])Y[\bar{v},-w^{-1}]\exp (wL[1]).
\end{equation*}%
Hence 
\begin{eqnarray*}
Z_{V }^{(2)}(\tau ,w,\rho ) &=&\sum_{m\geq 0}\sum_{v\in
V_{[m]}}\sum_{n\geq 0}\chi ^{n}\sum_{u\in V[n]}\langle E[\bar{v},w]\bar{u}%
,u\rangle Z_{V}^{(1)}(v,q) \\
&=&\sum_{m\geq 0}\sum_{v\in V_{[m]}}\Tr_{V}(E[\bar{v},w]\chi
^{L[0]})Z_{V}^{(1)}(v,q).
\end{eqnarray*}%
Now consider 
\begin{gather*}
\Tr_{V}(E[\bar{v},w]\chi ^{L[0]})= \\
w^{m}\sum\limits_{r,s\geq 0}(-1)^{r+m}\frac{1}{r!s!}\Tr_{V}(L[-1]^{r}\bar{v}%
[r-s-m-1]L[1]^{s}\chi ^{L[0]}).
\end{gather*}%
The leading term in $w$ is $w^{0}$ (arising from $\bar{v}=\vac $) and is given by 
\begin{equation*}
\Tr_{V}(E[\vac ,w]\chi ^{L[0]})=f(\chi )^{c/24}Z_{V}^{(1)}(f(\chi )).
\end{equation*}%
This follows from (\ref{Z1rhotrace}) and the isomorphism between the
original and square bracket formalisms. Taking $w\rightarrow 0$ for fixed $%
\chi $  the result follows.  $\hfill \square $

\section{The Heisenberg VOA}
\label{Sect_HVOA}
In this Section we compute the genus two partition function in the $\rho$-formalism
for the rank 1 Heisenberg VOA $M$. 
We also compute the genus two $n$-point function for $n$ copies of 
Heisenberg vector $a$ and the
genus two one point function for the Virasoro vector $\omega$. 
The main results mirror those obtained in the $\epsilon $-formalism in Section~6 of \cite{MT4}. 

\subsection{The genus two partition function $Z_{M}^{(2)}(\tau ,w,\rho )$}
\label{Subsect_Z2rho}
We begin by establishing a formula for
$Z_{M}^{(2)}(\tau ,w,\rho )$ in terms
of the infinite matrix $R$  \eqref{Rdef}.
\begin{theorem}
\label{Theorem_Z2_boson_rho} 
We have
\begin{equation}
Z_{M}^{(2)}(\tau ,w,\rho )=\frac{Z_{M}^{(1)}(\tau )}{\det (1-R)^{1/2}},
\label{Z2_1bos_rho}
\end{equation}
where $Z_{M}^{(1)}(\tau )=1/\eta (\tau )$.
\end{theorem}
\begin{remark}
\label{Rem_Z2mult} From Remark~2 of \cite{MT4} it follows that the
genus two partition function for $l$ free bosons $M^{l}$ is just the $l^{th}$
power of \eqref{Z2_1bos_rho}.
\end{remark}
\noindent
\textbf{Proof.} 
The proof is similar in structure to that of Theorem~5 of \cite{MT4}.
From \eqref{Z2_def_rho} we have
\begin{equation}
Z_{M}^{(2)}(\tau ,w,\rho )=
\sum_{n\geq 0}\sum_{u\in M_{[n]}}Z_{M}^{(1)}(u,\bar{u},w,\tau )\rho^{n},
\label{eq: part func II}
\end{equation}%
where $u$ ranges over a basis of $M_{[n]}$ and $\bar{u}$ is 
the dual state with respect to the square-bracket LiZ metric. 
$Z_{M}^{(1)}(u,v,w,\tau )$ is a genus one Heisenberg $2$-point function \eqref{Z1_2pt}.
We choose the square bracket Fock basis:
\begin{equation}
v=a[-1]^{e_{1}}\ldots   a[-p]^{e_{p}}\vac .  \label{eq: sq fock vec}
\end{equation}%
The Fock state $v$ naturally corresponds to an unrestricted partition 
$\lambda =\{1^{e_{1}}\ldots   p^{e_{p}}\}$ 
%with $\vert\lambda\vert=\sum_i e_{i}$ elements 
of $n=\sum_{1\leq i\leq p}ie_{i}$. 
We write $v=v(\lambda )$ to indicate this correspondence.
The Fock vectors form an orthogonal set from \eqref{eq: inner prod} with
\begin{equation*}
\bar v(\lambda )=\frac{1}{\prod_{1\leq i\leq p}(-i)^{e_{i}}e_{i}!} v(\lambda ).
\end{equation*}

\medskip The $2$-point function $Z_{M}^{(1)}(v(\lambda ),v(\lambda ),w,\tau) $ 
is given in Corollary~1 of \cite{MT1} where it is
denoted by $F_{M}(v,w_{1},v,w_{2};\tau )$. 
In order to describe this explicitly we introduce the set  
$\Phi _{\lambda,2 }$ which is the disjoint union of two isomorphic
label sets $\Phi_{\lambda }^{(1)}$, $\Phi _{\lambda }^{(2)}$ 
each with $e_{i}$ elements labelled $i$ determined by $\lambda$. 
Let $\iota:\Phi _{\lambda }^{(1)}\leftrightarrow \Phi _{\lambda }^{(2)}$ denote
the canonical label identification. 
Then we have
(loc. cit.) 
\begin{equation}
Z_{M}^{(1)}(v(\lambda ),v(\lambda ),w,\tau )=
Z_{M}^{(1)}(\tau )\sum_{\phi\in F(\Phi _{\lambda ,2})}\Gamma (\phi ),  \label{eq: genus1 2-point func}
\end{equation}%
where 
\begin{equation}
\Gamma (\phi ,w,\tau )=\Gamma (\phi )=\prod_{\{r,s\}}\xi (r,s,w,\tau ),
\label{eq: Gamma II}
\end{equation}%
and $\phi $ ranges over the elements of $F(\Phi _{\lambda ,2})$, the fixed-point-free involutions in $\Sigma (\Phi _{\lambda ,2})$ and where  $\{r,s\}$
ranges over the orbits of $\phi $ on $\Phi _{\lambda ,2}$. Finally 
\begin{equation}
\xi(r,s)=\xi (r,s,w,\tau )=\left\{ 
\begin{array}{ll}
C(r,s,\tau ), & \mbox{if}\ \{r,s\}\subseteq \Phi _{\lambda }^{(a)},a=1\ %
\mbox{or}\ 2, 
\notag
\\ 
D(r,s,w_{ab},\tau ) & \mbox{if}\ r\in \Phi _{\lambda }^{(a)},s\in \Phi
_{\lambda }^{(b)},a\neq b.%
\end{array}%
\right.,
\label{xirs}
\end{equation}
where  $w_{12}=w_{1}-w_{2}=w$ and $w_{21}=w_{2}-w_{1}=-w$.

\begin{remark}
\label{Rem_Dsym} Note that $\xi $ is well-defined since $D(r,s,w_{ab},\tau
)=D(s,r,w_{ba},\tau )$.
\end{remark}
 
Using the expression \eqref{eq: genus1 2-point func}, it follows that the genus two partition function
\eqref{eq: part func II} can be expressed as 
\begin{equation}
Z_{M}^{(2)}(\tau ,w,\rho )=Z_{M}^{(1)}(\tau )\sum\limits_{\lambda=\{ i^{e_{i}}\}} 
\frac{E(\lambda)}{\prod_{i}(-i)^{e_{i}}e_{i}!}
\rho^{\sum i e_{i}},
\label{ZMlambda}
\end{equation}
where $\lambda$ runs over all unrestricted partitions and
\begin{eqnarray}
E(\lambda)=\sum_{\phi\in F(\Phi _{\lambda ,2})}\Gamma (\phi ).
\label{Elambda}
\end{eqnarray}

\medskip 
We employ the doubly-indexed diagrams of Section~\ref{Subsect_Graphs}.
Consider the \lq canonical\rq\ matching defined by $\iota $ as a fixed-point-free involution.
We may then compose $\iota $ with each fixed-point-free involution 
$\phi\in F(\Phi _{\lambda ,2})$ to define a 1-1 mapping $\iota\phi $ on the underlying labelled set $\Phi _{\lambda ,2}$.
For each $\phi$ we define a doubly-indexed diagram  $D$ 
whose nodes are labelled by $k,a$ for an element 
$k\in \Phi_{\lambda }^{(a)}$ for $a=1,2$ and with
cycles corresponding to the orbits of the cyclic group $\langle \iota \phi\rangle $. 
Thus, if $l=\phi(k)$ for $k\in \Phi _{\lambda }^{(a)}$ and
 $l\in \Phi _{\lambda}^{(\bar{b})}$ and 
 $\iota :\bar{b}\mapsto b $ with convention \eqref{bar}
 then the corresponding doubly-indexed diagram contains the edge 
\begin{center}
\begin{picture}(250,15)
\put(100,5){\line(1,0){40}}
\put(84,10){\makebox(0,0){$(k,a)$}}
\put(100,5){\circle*{4}}
\put(157,10){\makebox(0,0){$(l,b)$}}
\put(140,5){\circle*{4}}
\end{picture}

{\small Fig.~3 A Doubly-Indexed Edge}
\end{center}
Consider the permutations of $\Phi_{\lambda ,2}$ that commute with $\iota $ and preserve
both $\Phi_{\lambda }^{(1)}$ and $\Phi_{\lambda }^{(2)}$. 
We denote this
group, which is plainly isomorphic to $\Sigma (\Phi_{\lambda })$, by $%
\Delta_{\lambda }$. By definition, an automorphism of a doubly-indexed diagram $%
D $ in the above sense is an element of $\Delta_{\lambda }$ which preserves
edges and node labels.

\medskip For a doubly-indexed diagram $D$ corresponding to the partition $\lambda
=\{1^{e_{1}}\ldots  p^{e_{p}}\}$ we set 
\begin{equation}
\gamma (D)=\frac{\prod_{\{k,l\}}\xi (k,l,w,\tau )}{\prod (-i)^{e_{i}}}\rho
^{\sum ie_{i}}  \label{eq: diagram weight II}
\end{equation}%
where $\{k,l\}$ ranges over the edges of $D$. We now have all the
pieces assembled to copy the arguments used to prove Theorem~5 of \cite{MT4}. 
First we find 
\begin{equation}
\sum\limits_{\lambda=\{ i^{e_{i}}\}} 
\frac{E(\lambda)}{\prod_{i}(-i)^{e_{i}}e_{i}!}
\rho^{\sum i e_{i}}=
\sum_{D}\frac{\gamma (D)}{|%
\mathrm{Aut}(D)|},  \label{eq: Z2Mrho}
\end{equation}%
the sum ranging over all doubly-indexed diagrams.

\medskip We next introduce a weight function $\zeta $ as follows: 
for a doubly-indexed diagram $D$ we set $\zeta (D)=\prod \zeta
(E)$, the product running over all edges. Moreover for an
edge $E$ with nodes labelled $(k,a)$ and $(l,b)$ as in Fig.~3, we set 
\begin{equation*}
\zeta (E)=R_{ab}(k,l),
\end{equation*}%
for $R$ of \eqref{Rdef}. We then find
\begin{lemma}
\label{Lemma_om_gamma_R}
$\zeta (D)=\gamma (D)$.
\end{lemma}
\noindent \textbf{Proof.}
From (\ref{eq: diagram weight II}) it follows that
for a doubly-indexed diagram $D$  we have 
\begin{equation}
\gamma (D)=\prod_{\{k,l\}}-\frac{\xi (k,l,w,\tau )\rho ^{(k+l)/2}}{\sqrt{kl}}
,  \label{eq: gammaprod II}
\end{equation}%
the product ranging over the edges $\{k,l\}$ of $D$. 
So to prove the Lemma it suffices to show that if
$k,l$ lie in $\Phi _{\lambda }^{(a)},\Phi _{\lambda }^{(\bar{b})}$
respectively then the $(a,b)$-entry of $R(k,l)$ coincides with the
corresponding factor of (\ref{eq: gammaprod II}). This follows from our
previous discussion together with Remark \ref{Rem_Dsym}. \quad  $\hfill \square $

\medskip From Lemma \ref{Lemma_om_gamma_R} and following similar arguments to 
the proof of Theorem~5 of \cite{MT4} we find 
\begin{eqnarray*}
\sum_{D}\frac{\gamma (D)}{|\mathrm{Aut}(D)|}
&=&\sum_{D}\frac{\zeta (D)}{|\mathrm{Aut}(D)|}\\
&=&\exp\left(\sum_{L}\frac{\zeta (L)}{|\mathrm{Aut}(L)|}\right),
\end{eqnarray*}%
where $L$ denotes the set of non-isomorphic unoriented doubly indexed cycles. 
Orient these cycles, say in a clockwise direction.
Let $\{M\}$ denote the set of non-isomorphic oriented doubly indexed cycles and $\{M_{n}\}$ the oriented cycles with $n$ nodes. Then we find (cf.\  \cite{MT4}, Lemma 2) that 
\begin{equation*}
\frac{1}{n}\mathrm{Tr} R^n=\sum_{M_{n}}\frac{\zeta (M_{n})}{|\mathrm{Aut}(M_{n})|}.
\end{equation*}
It follows that
\begin{eqnarray*}
\sum_{L}\frac{\zeta (L)}{|\mathrm{Aut}(L)|}
&=& 
\frac{1}{2}\sum_{M}\frac{\zeta (M)}{|\mathrm{Aut}(M)|}\\
&=& 
\frac{1}{2}\mathrm{Tr} \left(\sum\limits_{n\ge 1}\frac{1}{n}R^n\right)\\
&=& 
-\frac{1}{2}\mathrm{Tr} \log(I-R)\\
&=& 
-\frac{1}{2}\log\det(I-R).
\end{eqnarray*}
This completes the proof of Theorem \ref{Theorem_Z2_boson_rho}. $\hfill \Box$
\medskip

We may also find a product formula analogous to Theorem~6 of \cite{MT4}. Let $\mathcal{R}$ denote the rotationless doubly-indexed oriented cycles i.e. cycles with trivial automorphism group. Then we find 
\begin{theorem}
\label{Theorem_Z2_boson_prod_rho}
\begin{equation}
Z_{M}^{(2)}(\tau ,w,\rho ) = \frac{Z_{M}^{(1)}(\tau )}{\prod_{\mathcal{R}%
}(1-\zeta (N))^{1/2}} 
  \label{eq: prodIII}
\end{equation}
$\hfill \square $
\end{theorem}

\subsection{Holomorphic and modular-invariance properties}
\label{Subsect_HolMod}
In Section~\ref{Subsect_rho} we reviewed the genus two $\rho $-sewing formalism and
introduced the domain $\mathcal{D}^{\rho }$ which parametrizes the genus
two surface. An immediate consequence of Theorem \ref{Theorem_R} is the
following.
\begin{theorem}
\label{Theorem_Z2_boson_rho_hol}$Z_{M}^{(2)}(\tau ,w,\rho )$ is holomorphic
in $\mathcal{D}^{\rho }$.   
\end{theorem}
$\hfill \square $

\medskip
We next consider the invariance properties of the genus two partition
function with respect to the action of the $\mathcal{D}^{\rho }$-preserving
group $\Gamma_{1}$ reviewed in Section~\ref{Subsect_rho}. 
Let $\chi $ be the character of 
$SL(2,\mathbb{Z})$ defined by its action on $\eta (\tau )^{-2}$, i.e. 
\begin{equation}
\eta (\gamma \tau )^{-2}=\chi (\gamma )\eta (\tau )^{-2}(c\tau +d)^{-1},
\label{eq: eta}
\end{equation}%
where $\gamma =\left( 
\begin{array}{cc}
a & b \\ 
c & d%
\end{array}%
\right) \in SL(2,\mathbb{Z})$. Recall (e.g. \cite{Se}) that $\chi (\gamma )$
is a twelfth root of unity. For a function $f(\tau )$ on $\mathbb{H}%
_{1},k\in \mathbb{Z}$ and $\gamma \in SL(2,\mathbb{Z})$, we define 
\begin{equation}
f(\tau )|_{k}\gamma =f(\gamma \tau )\ (c\tau +d)^{-k},  \label{slashaction}
\end{equation}%
so that 
\begin{equation}
Z_{M^{2}}^{(1)}(\tau )|_{-1}\gamma =\chi (\gamma )Z_{M^{2}}^{(1)}(\tau ).  \label{Z1modgam}
\end{equation}
At genus two, analogously to (\ref{slashaction}), we define 
\begin{equation}
f(\tau ,w,\rho  )|_{k}\gamma =f(\gamma (\tau ,w,\rho ))
\det (C\Omega +D)^{-k}.  \label{eq: Gaction}
\end{equation}%
Here, the action of $\gamma $ on the right-hand-side is as in (\ref{Drho}).
We have abused notation by adopting the following conventions in (\ref{eq:
Gaction}), which we continue to use below: 
\begin{equation}
\Omega =F^{\rho }(\tau ,w,\rho),\ \gamma =\left( 
\begin{array}{cc}
A & B \\ 
C & D%
\end{array}%
\right) \in Sp(4,\mathbb{Z})  \label{Omegaconvention}
\end{equation}%
where $F^{\rho }$ is as in Theorem \ref{Theorem_period_rho}, and $\gamma 
$ is identified with an element of $Sp(4,\mathbb{Z})$ via 
\eqref{Gamma1} and \eqref{eq: modtrans}. Note that \eqref{eq: Gaction} defines a right action of $%
G $ on functions $f(\tau ,w,\rho)$. 
We then have a natural analog of Theorem~8 of \cite{MT4}
\begin{theorem}
\label{Theorem_Z2_rho_G1}
If $\gamma \in \Gamma _{1}$ then 
\begin{equation*}
Z_{M^{2}}^{(2)}(\tau ,w,\rho )|_{-1}\gamma =
\chi (\gamma)Z_{M^{2}}^{(2)}(\tau ,w,\rho ).
\end{equation*}
\end{theorem}
\begin{corollary}
\label{Cor_Z2_24_eps_G1}
If $\gamma \in \Gamma _{1}$ with $Z_{M^{24}}^{(2)}=(Z_{M^{2}}^{(2)})^{12}$ then 
\begin{equation*}
Z_{M^{24}}^{(2)}(\tau ,w,\rho )|_{-12}\gamma =Z_{M^{24}}^{(2)}(\tau ,w,\rho ).
\end{equation*}
\end{corollary}
\noindent \textbf{Proof.}
The proof is similar to that of Theorem~8 of \cite{MT4}.
We have to show that 
\begin{equation}
Z_{M^{2}}^{(2)}(\gamma .(\tau ,w,\rho ))\det (C\Omega +D)=\chi
(\gamma )Z_{M^{2}}^{(2)}(\tau ,w,\rho )
\label{eq: Z^2_rho_identity}
\end{equation}%
for $\gamma \in \Gamma _{1}$ where $\det (C\Omega _{11}+D)=c_{1}\Omega
_{11}+d_{1}$. 
Consider the determinant formula (\ref{Z2_1bos_rho}). For $\gamma \in \Gamma
_{1}$ define 
\begin{equation*}
R_{ab}^{\,\prime }(k,l,\tau ,w,\rho )=R_{ab}(k,l,\frac{a_{1}\tau +b_{1}}{%
c_{1}\tau +d_{1}},\frac{w}{c_{1}\tau +d_{1}},\frac{\rho }{(c_{1}\tau
+d_{1})^{2}})
\end{equation*}%
following (\ref{gam1_rho}). We find from Section 6.3 of \cite{MT2} that 
\begin{eqnarray*}
1-R^{\,\prime } &=&1-R-\kappa \Delta \\
&=&(1-\kappa S).(1-R),
\end{eqnarray*}%
where 
\begin{eqnarray*}
\Delta _{ab}(k,l) &=&\delta _{k1}\delta _{l1}, \\
\kappa &=&\frac{\rho }{2\pi i}\frac{c_{1}}{c_{1}\tau +d_{1}}, \\
S_{ab}(k,l) &=&\delta _{k1}\sum_{c\in \{1,2\}}((1-R)^{-1})_{cb}(1,l).
\end{eqnarray*}%
Since $\det (1-R)$ and $\det (1-R^{\,\prime })$ are convergent on $\mathcal{D}%
^{\rho }$ we find%
\begin{equation*}
\det (1-R^{\,\prime })=\det (1-\kappa S).\det (1-R).
\end{equation*}%
Indexing the columns and rows by $(a,k)=(1,1),(2,1),\ldots   (1,k),(2,k)\ldots   
$ and noting that $S_{1b}(k,l)=S_{2b}(k,l)$ we find that 
\begin{eqnarray*}
\det (1-\kappa S) &=&\left\vert 
\begin{array}{cccc}
1-\kappa S_{11}(1,1) & -\kappa S_{12}(1,1) & -\kappa S_{11}(1,2) & \cdots \\ 
-\kappa S_{11}(1,1) & 1-\kappa S_{12}(1,1) & -\kappa S_{11}(1,2) & \cdots \\ 
0 & 0 & 1 & \cdots \\ 
\vdots & \vdots & \vdots & \ddots%
\end{array}%
\right\vert \\
&=&1-\kappa S_{11}(1,1)-\kappa S_{12}(1,1), \\
&=&1-\kappa \sigma ((1-R)^{-1})(1,1)),
\end{eqnarray*}%
where $\sigma (M)$ denotes the sum of the entries of a (finite)
matrix $M$. Applying (\ref{Om11rho}), it is clear that 
\begin{equation*}
\det (1-\kappa S)=\frac{c_{1}\Omega _{11}+d_{1}}{c_{1}\tau +d_{1}}.
\end{equation*}%
The Theorem then follows from (\ref{Z1modgam}). $\hfill \square $

\begin{remark}
\label{Rem_Drhohat_mod} $Z_{M^{2}}^{(2)}(\tau ,w,\rho )$ can be
trivially considered as function on the covering space $\mathcal{\hat{D}}%
^{\rho }$ discussed in \cite{MT2}, Section 6.3.
Then $Z_{M^{2}}^{(2)}(\tau ,w,\rho )$ is modular with respect to $L=\hat{H}%
\Gamma _{1}$ with trivial invariance under the action of the Heisenberg
group $\hat{H}$ (loc.\ cit.).
\end{remark}

\subsection{Some genus two $n$-point functions}
\label{Subsect_HVOA_npt} 
In this section we calculate some examples of
genus two $n$-point functions for the rank one Heisenberg VOA $M$. 
We consider here the examples of the $n$-point function for the Heisenberg vector $a$
and the 1-point function for the Virasoro vector $\tilde\omega$. 
We find that, up to an overall factor of the partition function,  the formal differential form associated with the Heisenberg 
$n$-point function is described in terms of the global symmetric two form $\omega^{(2)}$ \cite{TUY} 
whereas the Virasoro 1-point function is described
by the genus two projective connection \cite{Gu}. 
These results agree with those found in \cite{MT4} in the $\epsilon$-formalism up to an overall $\epsilon$-formalism partition function factor.

\medskip
The genus two Heisenberg vector 1-point function with the Heisenberg vector $a$ inserted at $x$ is
 $Z_{M}^{(2)}(a,x;\tau,w,\rho)=0$ since 
$Z_{M}^{(1)}(Y[a,x]Y[v,w]v,\tau)=0$ from \cite{MT1}.
The 2-point function for two Heisenberg vectors inserted at
 $x_1,x_2$ is 
\begin{equation}
Z_{M}^{(2)}(a,x_1;a,x_2;\tau,w,\rho)=
\sum_{r\geq 0}\rho^{r}\sum_{v\in M_{[r]}}
Z_{M}^{(1)}(a,x_1;a,x_2;v,w_{1},\bar{v},w_{2};\tau).
%&=&
%\sum_{r\geq 0}\rho^{r}\sum_{v\in M_{[r]}}
%Z_{M}^{(1)}(Y[a,x_1]Y[a,x_2]Y[v,w]\bar{v};\tau),\notag\\
\label{Z2a1a1}
\end{equation}%
We consider the associated formal differential form
\begin{equation}
\mathcal{F}_{M}^{(2)}(a,a;\tau,w,\rho)=Z_{M}^{(2)}(a,x_1;a,x_2;\tau,w,\rho)dx_{1}dx_{2},
\label{eq:F2aa}
\end{equation}
and find that it is determined by the bilinear
form $\omega^{(2)}$ (\ref{omegag}):
\begin{theorem}
\label{theorem:Z2aa} The genus two Heisenberg vector 2-point function is given by
\begin{eqnarray}
\mathcal{F}_{M}^{(2)}(a,a;\tau,w,\rho)=
\omega^{(2)}(x_{1},x_{2})Z_{M}^{(2)}(\tau,w,\rho).  \label{eq: F2aa}
\end{eqnarray}
\end{theorem}
\noindent \textbf{Proof.} 
The proof proceeds along the same lines as Theorem~\ref{Theorem_Z2_boson_rho}.
As before, we let $v(\lambda )$ denote a Heisenberg
Fock vector (\ref{eq: sq fock vec}) determined by an unrestricted partition $%
\lambda =\{1^{e_{1}}\ldots   p^{e_{p}}\}$ with label set $\Phi_{\lambda}$.
Define a label set for the four vectors $a,a,v(\lambda ),v(\lambda )$ given by $%
\Phi=\Phi_1\cup \Phi_2\cup \Phi_{\lambda}^{(1)} \cup \Phi_{\lambda}^{(2)}$ 
for $\Phi_1,\Phi_2=\{ 1\}$ and 
let $F(\Phi)$ denote the set of fixed point free involutions on $\Phi$. 
For $\phi=\ldots  (r s)\ldots  \in F(\Phi)$ let 
$\Gamma(x_1,x_2,\phi)=\prod_{(r,s)}\xi(r,s)$  as defined in 
\eqref{xirs} for $r,s\in\Phi_{\lambda}^{(2)}=\Phi_{\lambda}^{(1)}\cup \Phi_{\lambda}^{(2)}$ and 
\begin{equation}
\xi (r,s)=\left \{ 
\begin{array}{ll}
D(1,1,x_i-x_j,\tau )=P_2(\tau,x_i-x_j), &  r,s\in \Phi_{i},\ i\neq j,
\\ 
D(1,s,x_{i}-w_{a},\tau )=sP_{s+1}(\tau,x_{i}-w_{a}), & r\in \Phi_{i},\ s\in \Phi_{\lambda}^{(a)},
%\\ 
%D(r,1,w_{a}-x_{i},\tau )=rP_{r+1}(\tau,w_{a}-x_{i}), & r\in \Phi_{\lambda}^{(a)};\ s\in \Phi_{i},
\end{array}
\right.  \label{eq: X12 weight}
\end{equation}
for $i,j,a\in\{1,2\}$ with $D$ of \eqref{Dkldef}. 
Then following Corollary~1
of \cite{MT1} we have 
\begin{equation*}
Z_{M}^{(1)}(a,x_1; a,x_2; v(\lambda ),w_{1};v(\lambda ),w_{2}; \tau )=Z_{M}^{(1)}(\tau )
\sum_{\phi \in F(\Phi)}\Gamma (x_1,x_2,\phi ).
\end{equation*}%
We then obtain the following analog of (\ref{ZMlambda}) 
\begin{equation}
\mathcal{F}^{(2)}(a,a;\tau,w,\rho)=
Z_{M}^{(1)}(\tau)
\sum_{\lambda =\{i^{e_{i}}\}}%
\frac{E(x_1,x_2,\lambda)}{\prod_{i}(-i)^{e_{i}}e_{i}!}
\rho^{\sum ie_{i}}dx_1dx_2,  \label{eq: part func lambda a1a2}
\end{equation}%
where 
\begin{equation*}
E(x_1,x_2,\lambda ) = \sum_{\phi\in F(\Phi)}\Gamma(x_1,x_2,\phi).
\end{equation*}

The sum in \eqref{eq: part func lambda a1a2} can be re-expressed as the sum
of weights $\zeta (D)$ for isomorphism classes of doubly-indexed
configurations $D$ where here $D$ includes
two distinguished valency one nodes labelled $x_i$ (see Section~\ref{Subsect_Graphs}) 
corresponding to the label sets $\Phi_1,\Phi_2=\{ 1\}$. 
As before, $\zeta (D)=\prod_E\zeta(E)$ for standard doubly-indexed
 edges $E$ augmented by the contributions from edges connected to the
two valency one nodes with weights as in \eqref{eq: zeta1a2}.
Thus we find
\begin{equation*}
\mathcal{F}^{(2)}(a,a;\tau,w,\rho)=
Z_{M}^{(1)}(\tau)
\sum_{D}\frac{\zeta (D)}{\prod_{i}e_{i}!}dx_1dx_2,
\end{equation*}
Each $D$ can be decomposed into \emph{exactly} one necklace configuration $N$ of type $\mathcal{N}(x;y)$ of \eqref{Nxydef} 
connecting the two distinguished nodes and a standard configuration ${\hat D}$ of the type
appearing in the proof of Theorem~\ref{Theorem_Z2_boson_rho} with $\zeta (D)=\zeta(N)\zeta ({\hat D})$. Since $|\mathrm{Aut}(N)|=1$ we obtain 
\begin{eqnarray*}
\mathcal{F}^{(2)}(a,a;\tau,w,\rho)&=&
Z_{M}^{(1)}(\tau)
\sum_{{\hat D}}\frac{\zeta ({\hat D})}{|\mathrm{Aut}({\hat D})|} 
\sum_{N\in \mathcal{N}(x;y)}\zeta(N)\\
&=&Z_{M}^{(2)}(\tau,w,\rho )\zeta(x_1;x_2) \\
&=&Z_{M}^{(2)}(\tau,w,\rho )\omega^{(2)}(x_1,x_2),
\end{eqnarray*}
using (\ref{om_graph}) of Proposition \ref{Propomeganugraph}. 
 $\hfill \square $

\medskip
Theorem \ref{theorem:Z2aa} result can be generalized to compute the $n$-point function corresponding
to the insertion of
$n$ Heisenberg vectors. We find that it
vanishes for $n$ odd, and for $n$ even is determined
by the symmetric tensor 
\begin{equation}
\mathrm{Sym}_n\omega^{(2)}= \sum_{\psi}\prod_{(r,s)}\omega^{(2)}(x_r,x_s),
\label{eq: bos n form}
\end{equation}
where the sum is taken over the set of fixed point free involutions 
$\psi=\ldots   (rs)\ldots  $ of the labels $\{1,\ldots  ,n\}$. We then have
\begin{theorem}
\label{theorem:Z2an} The genus two Heisenberg vector $n$-point function vanishes for odd 
$n$ even; for even $n$ it is
given by the global symmetric meromorphic $n$-form:
\begin{equation}
\frac{\mathcal{F}_{M}^{(2)}(a,\ldots  ,a;\tau,w,\rho)}
{ Z_{M}^{(2)}(\tau,w,\rho)}
=\mathrm{Sym}_n\omega^{(2)}.  \label{eq: F2n}
\end{equation}
$\hfill \Box$
\end{theorem}
This agrees with the corresponding ratio in Theorem~10 of \cite{MT4} in the $\epsilon$-formalism 
and earlier results in \cite{TUY} based on an assumed analytic structure for the $n$-point function.

\medskip
Using this result and the associativity of vertex operators,  we can compute all $n$-point functions.
In particular, the 1-point function for the Virasoro vector 
${\tilde \omega}=\frac{1}{2}a[-1]a$ is as follows (cf.\ \cite{MT4}, Proposition 8):
\begin{proposition}
\label{Prop:Z2omega} 
The genus two 1-point function for the Virasoro vector $\tilde\omega$ inserted at $x$ is given by 
\begin{eqnarray}
\frac{\mathcal{F}_{M}^{(2)}({\tilde\omega};\tau,w,\rho)}
{Z_{M}^{(2)}(\tau,w,\rho)}
= \frac{1}{12}s^{(2)}(x),  \label{eq: F2omega}
\end{eqnarray}
where 
$s^{(2)}(x)=6\lim_{x \rightarrow y} \left (\omega^{(2)}(x,y)-\frac{dxdy}{%
(x-y)^2}\right )$
is  the genus two projective connection (\cite{Gu}). $\hfill \square$
\end{proposition}

\section{Lattice VOAs}
\label{Sect_LVOA}

\subsection{The genus two partition function $Z_{V_{L}}^{(2)}(\tau ,w,\rho )$}
\label{Subsect_Z2Lrho}
Let $L$ be an even lattice with $V_{L}$ the corresponding lattice theory
vertex operator algebra. The underlying Fock space is 
\begin{equation}
V_{L}=M^{l}\otimes C[L]=\oplus _{\beta \in L}M^{l}\otimes e^{\beta },
\label{lattice fock space}
\end{equation}%
where $M^{l}$ is the corresponding Heisenberg free boson theory of rank $l=$
dim $L$ based on $H=C\otimes _{Z}L$. We follow Section~\ref{Subsect_VOA} and \cite{MT1}
concerning further notation for lattice theories.
We utilize the Fock space $\{u\otimes e^{\beta }\}$ where $\beta $ ranges over $L$ and $u$ ranges
over the usual orthogonal basis for $M^{l}$. From Lemma \ref{Lemma_LiZ_lattice} and Corollary \ref{Corollary_cocycle_choice} we see that 
\begin{eqnarray}
Z_{V_{L}}^{(2)}(\tau ,w,\rho )&=&
\sum_{\alpha ,\beta \in L}
Z_{\alpha ,\beta}^{(2)}(\tau ,w,\rho ),  \label{eq: rk1rholattpartfunc}
 \\
Z_{\alpha ,\beta}^{(2)}(\tau ,w,\rho )&=&\sum_{n\geq 0}
\sum_{u\in M_{[n]}^{l}}Z_{M^{l}\otimes e^{\alpha }}^{(1)}(u\otimes e^{\beta },\overline{u}\otimes e^{-\beta
},w,\tau )\rho ^{n+(\beta ,\beta )/2}.
\label{eq: Zalbet}
\end{eqnarray}%
The general shape of the $2$-point function occurring in (\ref{eq: Zalbet}) is discussed extensively in \cite{MT1}. By Proposition
1 (loc. cit.) it splits as a product 
\begin{gather}
Z_{M^{l}\otimes e^{\alpha }}^{(1)}(u\otimes e^{\beta },u\otimes e^{-\beta
},w,\tau )=  \notag \\
Q_{M^{l}\otimes e^{\alpha }}^{\beta }(u,u,w,\tau )Z_{M^{l}\otimes e^{\alpha
}}^{(1)}(e^{\beta },e^{-\beta },w,\tau ),  \label{eq: rhosplitform}
\end{gather}%
where we have identified $e^{\beta }$ with $\vac \otimes e^{\beta }$,
and where $Q_{M^{l}\otimes e^{\alpha }}^{\beta }$ is a function\footnote{%
Note: in \cite{MT1} the functional dependence on $\beta $, here denoted by
a superscript, was omitted.} that we will shortly discuss in greater detail.
In \cite{MT1}, Corollary~5 (cf.\ the Appendix to the present paper) we
established also that 
\begin{equation}
Z_{M^{l}\otimes e^{\alpha }}^{(1)}(e^{\beta },e^{-\beta },w,\tau )=\epsilon
(\beta ,-\beta )q^{(\alpha ,\alpha )/2}\frac{\exp ((\alpha ,\beta )w)}{%
K(w,\tau )^{(\beta ,\beta )}}Z_{M^{l}}^{(1)}(\tau ),
\label{eq: latt2ptform}
\end{equation}%
where, as usual, we are taking $w$ in place of $z_{12}=z_{1}-z_{2}$. With
cocycle choice $\epsilon (\beta ,-\beta )=(-1)^{(\beta ,\beta )/2}$ (cf.\ Appendix) we may then rewrite (\ref{eq: Zalbet}) as 
\begin{eqnarray}
Z_{\alpha,\beta}^{(2)}(\tau ,w,\rho )&=& Z_{M^{l}}^{(1)}(\tau )
\exp \left\{\pi i\Big((\alpha ,\alpha)\tau 
+2(\alpha ,\beta )\frac{w}{2\pi i}
+\frac{(\beta ,\beta )}{2\pi i}%
\log \left(\frac{-\rho }{K(w,\tau )^{2}}\right)\Big)\right\}  \notag \\
&&\cdot \sum_{n\geq 0}\sum_{u\in M_{[n]}^{l}}
Q_{M^{l}\otimes e^{\alpha }}^{\beta }(u,\bar{u},w,\tau )
\rho ^{n}.  \label{eq: newrholattpartfunc}
\end{eqnarray}%
We note that this expression is, as it should be, independent of the choice
of branch for the logarithm function. We are going to establish the \emph{%
precise} analog of Theorem~14 of \cite{MT4} as follows:
\begin{theorem}
\label{Theorem_Z2_L_rho}
We have 
\begin{equation}
Z_{V_{L}}^{(2)}(\tau ,w,\rho )=Z_{M^{l}}^{(2)}(\tau ,w,\rho )\theta
_{L}^{(2)}(\Omega ),
\label{ZL2}
\end{equation}
where $\theta _{L}^{(2)}(\Omega )$ is the genus two theta function of $L$ (\cite{Fr}).
\end{theorem}
{\bf Proof.} We note that
\begin{equation}
\theta _{L}^{(2)}(\Omega )=\sum_{\alpha ,\beta \in L}\exp (\pi i((\alpha ,\alpha)\Omega _{11}+2(\alpha ,\beta )\Omega _{12}+(\beta
,\beta  )\Omega
_{22})).  \label{eq: ThetaL}
\end{equation}

\medskip
We first handle the case of rank $1$ lattices and then consider the general case. 
The inner double sum in \eqref{eq: newrholattpartfunc} is the object which requires attention, and we can
begin to deal with it along the lines of previous Sections. Namely,
arguments that we have already used several times show that the double sum
may be written in the form 
\begin{equation*}
\sum_{D}\frac{\gamma (D)}{|\mathrm{Aut}(D)|}=\exp \left( \frac{1}{2}%
\sum_{N\in \mathcal{N}}\gamma (N)\right) .
\end{equation*}%
Here, $D$ ranges over the oriented doubly indexed cycles of Section~\ref{Sect_HVOA},
while $N$ ranges over oriented doubly-indexed necklaces  
${\cal N}=\{{\cal N}(k,a;l,b)\}$ of \eqref{Nkl}. 
Leaving aside the definition of $\gamma (N)$ for now, we
recognize as before that the piece involving only connected diagrams with no
end nodes splits off as a factor.
Apart from a $Z_{M}^{(1)}(\tau )$ term this factor is, of course, precisely
the expression (\ref{eq: Z2Mrho}) for $M$. With these
observations, we see from (\ref{eq: newrholattpartfunc}) that the following
holds: 
\begin{eqnarray}
\frac{Z_{\alpha,\beta}^{(2)}(\tau ,w,\rho )}{Z_{M}^{(2)}(\tau ,w,\rho)}
&=&
\exp \left\{ 
i\pi\left((\alpha ,\alpha )\tau +2(\alpha,\beta )\frac{w}{2\pi i}+\frac{(\beta ,\beta )}{2\pi i}
\log \left(\frac{-\rho }{K(w,\tau )^{2}}\right)
\right.\right.  \notag 
\\
&&
\left.\left.
\quad+\frac{1}{2\pi i}\sum_{N\in \cal{N}}\gamma (N)\right)\right\}.  \label{eq: lattrhopartqnt}
\end{eqnarray}%
To prove Theorem~\ref{Theorem_Z2_L_rho}, we see from (\ref{eq: ThetaL}) and
\eqref{eq: lattrhopartqnt} that it is sufficient to establish that for each
pair of lattice elements $\alpha ,\beta \in L$, we have 
\begin{eqnarray}
\frac{Z_{\alpha,\beta}^{(2)}(\tau ,w,\rho )}{Z_{M}^{(2)}(\tau ,w,\rho)}=\exp\left(
\pi i\left((\alpha ,\alpha )\Omega _{11}
+2(\alpha ,\beta )\Omega_{12}
+(\beta ,\beta )\Omega _{22}\right)
\right).  \label{eq: rhoomegas}
\end{eqnarray}

Recall the formula for  $\Omega$ in Proposition~\ref{Prop_rhoperiod_graph}.
In order to reconcile (\ref{eq: rhoomegas}) with the formula for $\Omega$, we
must carefully consider the expression $\sum_{N\in \cal{N}}\gamma(N)$. 
The function $\gamma $ is essentially (\ref{eq: diagram weight II}), except that we also
get contributions from the end nodes which are now present. Suppose that an
end node has label $k\in \Phi ^{(a)},a\in \{1,2\}$. 
Then according to Proposition~1 and display~(45) of \cite{MT1} (cf. \eqref{A.2} of the Appendix to the
present paper), the contribution of the end node is equal to 
\begin{gather}
\xi_{\alpha ,\beta }(k,a)=
\xi_{\alpha ,\beta }(k,a,\tau ,w,\rho )=  \notag \\
\left\{ 
\begin{array}{ll}
\frac{\rho ^{k/2}}{\sqrt{k}}
\left(a,\delta _{k1}\alpha +C(k,0,\tau )\beta
+D(k,0,w,\tau )(-\beta) \right), & a=1 \\ 
\frac{\rho ^{k/2}}{\sqrt{k}}
\left(a,\delta _{k1}\alpha +C(k,0,\tau )(-\beta)
+D(k,0,-w,\tau )\beta \right), & a=2%
\end{array}%
\right.  \label{eq: xirho}
\end{gather}%
together with a contribution arising from the $-1$ in the denominator of (%
\ref{eq: diagram weight II}) (we will come back to this point later). Using
(cf.\ \cite{MT1}, displays (6), (11) and (12)) 
\begin{eqnarray*}
D(k,0,-w,\tau ) &=&(-1)^{k+1}P_{k}(-w,\tau )=-P_{k}(w,\tau ), \\
C(k,0,\tau ) &=&(-1)^{k+1}E_{k}(\tau ),
\end{eqnarray*}%
we can combine the two possibilities in (\ref{eq: xirho}) as follows (recalling
that $E_{k}=0$ for odd $k$): 
\begin{equation}
\xi _{\alpha ,\beta }(k,a)= 
(a,\alpha)\rho^{1/2}\delta _{k1}
+(a,\beta)d_{\bar{a}}(k),
\label{eq: xialphabetaform}
\end{equation}
where $d_{a}(k)$ is given by \eqref{dk}.
We may then compute the weight for an oriented doubly-indexed necklace $N\in{\cal N}(k,a;l,b)$ 
\eqref{Nkl}.
Let $N^{\prime }$ denote the oriented necklace from which the
two end nodes and edges have been \emph{removed} (we refer to these as \emph{%
shortened} necklaces).
From (\ref{eq: xialphabetaform}) we see that the total
contribution to $\gamma (N)$ is 
\begin{eqnarray}
-\xi_{\alpha ,\beta }(k,a )\xi_{\alpha ,\beta }(l,b )\gamma (N^{\prime }) 
&=&-\left[
(\alpha,\alpha)\rho\delta_{k1}\delta_{l1}
+(\beta,\beta)d_{\bar{a}}(k)d_{\bar{b}}(l)\right.
\notag\\
&&
\left.
+(\alpha,\beta)\rho^{1/2}\left(d_{\bar{a}}(k)\delta _{l,1}+d_{\bar{b}}(l)\delta _{k,1}\right)
\right]\gamma (N^{\prime }),\notag\\ 
\label{eq: gammaNalphabetaform}
\end{eqnarray}
where we note that a sign $-1$ arises from each \emph{pair} of nodes, as
follows from (\ref{eq: diagram weight II}).  

\medskip
We next consider the terms in \eqref{eq: gammaNalphabetaform} corresponding to 
$(\alpha ,\alpha),(\alpha ,\beta )$ and $(\beta ,\beta )$ separately, and show that they
are precisely the corresponding terms on each side of (\ref{eq: rhoomegas}).
This will complete the proof of Theorem \ref{Theorem_Z2_L_rho} in the case
of rank $1$ lattices. 
From (\ref{eq: gammaNalphabetaform}), an $(\alpha ,\alpha)$ term arises only if the end node weights $k,l$ are both equal to $1$. 
Hence $\sum\gamma(N^{\prime})=\zeta(1;1)$ (cf.\ \eqref{om_rho_weights}),
where the sum ranges over shortened necklaces with end nodes of weight 
$1\in \Phi ^{(a)}$ and $1\in \Phi ^{(b)}$.
Thus 
using Proposition~\ref{Prop_rhoperiod_graph}, the total contribution to the right-hand-side
 of \eqref{eq: rhoomegas} is equal to 
\begin{equation}
\tpi\tau -\rho \zeta(1;1)=\tpi\Omega _{11}.
\label{Omega11seen}
\end{equation}%

\medskip Next, from (\ref{eq: gammaNalphabetaform}) we see that an $(\alpha
,\beta )$-contribution arises whenever at least one of the end nodes has
label $1$. If the labels of the end nodes are unequal then the shortened
necklace with the \emph{opposite} orientation makes an equal contribution.
The upshot is that we may assume that the end node to the right of the
shortened necklace has label $l=1\in \Phi ^{(\bar{b})}$, as long as we count
accordingly. 
We thus find $\sum\gamma(N^{\prime})=\zeta(d;1)$ (cf.\ \eqref{om_rho_weights}),
where the sum ranges over shortened necklaces with end nodes of weight 
$k\in \Phi ^{(a)}$ and $1\in \Phi ^{(b)}$.
Then Proposition~\ref{Prop_rhoperiod_graph} implies that the total contribution to
the $(\alpha ,\beta )$ term on the right-hand-side of (\ref{eq: rhoomegas}) is
\begin{equation*}
2w-2\rho^{1/2}\zeta(d;1)=2\Omega_{12},
\end{equation*}%
as required.

\medskip It remains to deal with the $(\beta ,\beta )$ term, the details
of which are very much along the lines as the case $(\alpha ,\beta )$ just
handled. A similar argument shows that the contribution to the $(\beta
,\beta )$-term from (\ref{eq: gammaNalphabetaform}) is equal to the expression
$-\zeta(d;\bar{d})$ of \eqref{om_rho_weights}. 
Thus the total contribution to the $(\beta ,\beta )$ term on the
right-hand-side of (\ref{eq: rhoomegas}) is 
\begin{equation*}
\log \left(\frac{-\rho }{K(w,\tau )^{2}}\right)-\zeta(d;\bar{d})=\Omega _{22},
\end{equation*}%
as in (\ref{Om22rho}). 
This  completes the proof of Theorem \ref{Theorem_Z2_L_rho} in the rank $1$ case.

\medskip
As for the general case
- we adopt the mercy rule and omit details! The reader who has progressed
this far will have no difficulty in dealing with the general case, which
follows by generalizing the calculations in the rank $1$ case
just considered. $\hfill \square $

\medskip The analytic and automorphic properties of $Z_{V_{L}}^{(2)}(\tau,w,\rho )$ can be deduced from Theorem \ref{Theorem_Z2_L_rho} using the known behaviour of  $\theta^{(2)}_L(\Omega)$ and the analagous results for $Z^{(2)}_{M^l}(\tau, w, \rho)$ established in
Section 5. We simply record
\begin{theorem}
\label{Theorem_Z2_lattice_rho_hol}
$Z_{V_{L}}^{(2)}(\tau ,w,\rho )$ is
holomorphic on\ the domain $\mathcal{D}^{\rho }$. $ \hfill \square $
\end{theorem}

\subsection{Some genus two $n$-point functions}
\label{Subsect_LVOA_npt}
In this section we consider the genus two $n$-point functions for $n$
Heisenberg vectors and the 1-point function for the Virasoro vector $%
\tilde\omega$ for a rank $l$ lattice VOA. 
The results are similar to those of Section~\ref{Subsect_HVOA_npt}
so that detailed proofs will not be given.

\medskip
Consider the $1$-point function for a Heisenberg vector $a_{i}$ inserted at $x$. We define the differential 1-form 
\begin{equation}
\mathcal{F}_{\alpha,\beta}^{(2)}(a_{i};\tau,w,\rho)=
\sum_{n\geq 0}\sum_{u\in M_{[n]}^{l}}
Z_{M^{l}\otimes e^{\alpha }}^{(1)}
(a_{i},x;u\otimes e^{\beta },w;\overline{u}\otimes e^{-\beta},0,\tau )\rho ^{n+(\beta ,\beta )/2}dx.
\label{eq:Fab}
\end{equation}
This can be expressed in terms of the genus two holomorphic 1-forms $\nu_{1},\nu_{2}$ of \eqref{nui_rho} in a similar way to Theorem~12 of \cite{MT4}. Defining $$\nu_{{i},\alpha,\beta}(x)=(a_{i},\alpha)\nu_{1}(x)+(a_{i},\beta)\nu_{2}(x),$$ 
we
find 
\begin{theorem}
\label{theorem:Z2a_alpha} 
\begin{eqnarray}
\mathcal{F}_{\alpha,\beta}^{(2)}(a_{i};\tau,w,\rho)
=\nu_{i,\alpha,\beta}(x)Z_{\alpha,\beta}^{(2)}(\tau,w,\rho).
  \label{eq: Z2a_alpha}
\end{eqnarray}
\end{theorem}
\textbf{Proof.} The proof proceeds along the same lines as Theorems~\ref
{theorem:Z2aa} and \ref{Theorem_Z2_L_rho} and Theorem~12 of (op.\ cit.). 
We find that 
\begin{equation*}
\mathcal{F}_{\alpha,\beta}^{(2)}(a_{i};\tau,w,\rho)=
Z_{M}^{(1)}(\tau)\sum_{D}\frac{\gamma (D)}{|\mathrm{Aut}(D)|}dx,
\end{equation*}
where the sum is taken over isomorphism classes of doubly-indexed configurations $D$ where, in this case, each configuration includes one distinguished
valence one node labelled by $x$ as in \eqref{eq: zeta1a2}. 
Each $D$ can be decomposed into exactly
one necklace configuration of type $\mathcal{N}(x;k,a)$ of 
(\ref{Nx1def}), standard configurations of the type appearing in 
Theorem~\ref{theorem:Z2aa} and necklace contributions as in Theorem~\ref{ZL2}. The result then follows
on applying \eqref{eq: xialphabetaform} and the graphical expansion for $\nu_1(x),\nu_2(x)$ of \eqref{nu1_graph} and \eqref{nu2_graph}. $\hfill \square$

\medskip
Summing over all lattice vectors, we find that the Heisenberg 1-point function vanishes for $V_{L}$. 
Similarly, one can generalize Theorem \ref{theorem:Z2an}
concerning the $n$-point function for $n$ Heisenberg vectors $a_{i_{1}},\ldots,a_{i_{n}}$. Defining
\begin{eqnarray*}
\mathcal{F}^{(2)}_{\alpha,\beta}(a_{i_{1}},\ldots,a_{i_{n}};\tau,w,\rho)
&=&\prod_{t=1}^{n}dx_{t}
\sum_{n\geq 0}\sum_{u\in M_{[n]}^{l}}\rho ^{n+(\beta ,\beta )/2}\cdot 
\\&&
Z_{M^{l}\otimes e^{\alpha }}^{(1)}
(a_{i_{1}},x_{1};\ldots;a_{i_{n}},x_{n};
u\otimes e^{\beta },w;
\overline{u}\otimes e^{-\beta},0,\tau ),
\end{eqnarray*} 
we obtain the analogue of Theorem~13 of (op.\ cit.):
\begin{theorem}
\label{theorem:Z2anModule}  
\begin{equation}
\mathcal{F}^{(2)}_{\alpha,\beta}(a_{i_{1}},\ldots,a_{i_{n}};\tau,w,\rho)
=
\mathrm{Sym}_n\big(\omega^{(2)},\nu_{i_{t},\alpha,\beta}\big)Z_{\alpha,\beta}^{(2)}(\tau,w,\rho),
\label{eq: F2n alpha}
\end{equation}
the symmetric product of $\omega^{(2)}(x_r,x_s)$ and $\nu_{i_{t},\alpha,\beta}(x_{t})$   defined by 
\begin{equation}
\mathrm{Sym}_n\big(\omega^{(2)},\nu_{i_{t},\alpha,\beta}\big)=
\sum_{\psi}\prod_{(r,s)}\omega^{(2)}(x_s,x_s)\prod_{(t)}\nu_{i_{t},\alpha,\beta}(x_t),
\label{eq: bos n form alpha}
\end{equation}
where the sum is taken over the set of involutions $\psi=\ldots (ij)\ldots
(k)\ldots$ of the labels $\{1,\ldots,n\}$. $\hfill \Box$
\end{theorem}
 
\medskip
We may also compute the genus two 1-point function for the Virasoro vector $\tilde\omega=\frac{1}{2}\sum_{i=1}^{l} a_i[-1]a_{i}$  using associativity of vertex operators as in
Proposition~\ref{Prop:Z2omega}. 
We  find that for a rank $l$ lattice,
\begin{eqnarray*}
\mathcal{F}_{\alpha,\beta}^{(2)}(\widetilde \omega;\tau,w,\rho) 
&\equiv&
\sum_{n\geq 0}\sum_{u\in M_{[n]}^{l}}
Z_{M^{l}\otimes e^{\alpha }}^{(1)}
(\widetilde \omega,x;u\otimes e^{\beta },w;\overline{u}\otimes e^{-\beta},0,\tau )\rho ^{n+(\beta ,\beta )/2}dx^2
\\
&=&
\left( \frac{1}{2}\sum_{i}\nu_{i,\alpha,\beta}(x)^2+\frac{l}{2}s^{(2)}(x) \right)
Z_{\alpha,\beta}^{(2)}(\tau,w,\rho)
\\
&=&
Z_{M^l}^{(2)}(\tau,w,\rho) 
\left(\mathcal{D} + \frac{l}{12}s^{(2)} \right )
e^{
i\pi \left((\alpha ,\alpha )\Omega _{11}
+2(\alpha ,\beta )\Omega_{12}
+(\beta ,\beta )\Omega _{22}\right)}.
\end{eqnarray*} 
Here, we used \eqref{eq: rhoomegas} and the differential operator (\cite{Fa}, \cite{U}, \cite{MT4}) 
\begin{equation}
\mathcal{D}=\frac{1}{2\pi i}\sum_{1\le a\le b\le 2}\nu_a\nu_b\frac{\partial}
{\partial \Omega_{ab}}.  \label{eq: Del 2form}
\end{equation}
Defining the normalized Virasoro $1$-point form
\begin{equation}
\hat{\mathcal{F}}_{V_{L}}^{(2)}(\tilde\omega;\tau,w,\rho)=
\frac{\mathcal{F}_{V_{L}}^{(2)}(\tilde\omega;\tau,w,\rho)}
{Z_{M^{l}}^{(2)}(\tau,w,\rho)},  \label{Zomega_hat_rho}
\end{equation}
we obtain 
\begin{proposition}
\label{Prop: Lattice Virasoro} The normalized Virasoro 1-point function for the 
lattice theory $V_L$
 satisfies  
\begin{equation}
\hat{\mathcal{F}}_{V_{L}}^{(2)}({\tilde\omega};\tau,w,\rho )
=
 \left(\mathcal{D} + \frac{l}{12}s^{(2)} %
\right )\theta _{L}^{(2)}(\Omega ).  \label{eq: F2omega lattice Del }
\end{equation}
$\hfill \Box$
\end{proposition}
(The Ward identity \eqref{eq: F2omega lattice Del } is similar to Proposition~11 in \cite{MT4} in the $\epsilon$-sewing formalism.)

\medskip
Finally, we can obtain the analogue of Proposition~12 (op.\ cit.), where we find that $\hat{\mathcal{F}}_{V_{L}}^{(2)}$
enjoys the same modular properties as $\hat{Z}_{V_{L}}^{(2)}=\theta _{L}^{(2)}(\Omega(\tau,w,\rho) )$.
That is, 
\begin{proposition}
\label{prop: F2omega lattice mod} The normalized Virasoro 1-point
function for a lattice VOA obeys
\begin{equation}
\hat{\mathcal{F}}_{V_{L}}^{(2)}
({\tilde\omega};\tau,w,\rho)|_{l/2}\gamma= 
\big(\mathcal{D} + \frac{l}{12}s^{(2)} \big )
\big(\hat{Z}_{V_{L}}^{(2)}(\tau,w,\rho )|_{l/2}\gamma\big),
\label{eq: F2omega lattice Del mod}
\end{equation}
for $\gamma \in \Gamma_{1}$. $\hfill \Box$
\end{proposition}

\section{Comparison between the $\epsilon$ and $\rho$-Formalisms}
\label{Sect_Comp}
In this section we consider the relationship between the genus two boson and
lattice partition functions computed in the $\epsilon $-formalism of \cite{MT4} (based on a sewing construction with two separate tori with modular parameters $\tau_{1},\tau_{2}$ and a sewing parameter $\epsilon$) and the $\rho$-formalism developed in this paper. 
We write 
\begin{eqnarray*}
Z_{V,\epsilon}^{(2)}
=Z_{V,\epsilon}^{(2)}(\tau _{1},\tau _{2},\epsilon )
&=&
\sum_{n\geq 0}\epsilon ^{n}
\sum_{u\in V_{[n]}}Z_{V}^{(1)}(u,\tau_{1})
 Z_{V}^{(1)}(\bar{u},\tau_{2} ),\\
Z_{V,\rho }^{(2)}
=Z_{V,\rho }^{(2)}(\tau ,w,\rho )
&=&
\sum_{n\geq 0}\rho ^{n}\sum_{u\in
V_{[n]}}Z_{V}^{(1)}(\bar{u},u,w,\tau ).
\end{eqnarray*} 
Although, for a given VOA $V$, the
partition functions enjoy many similar properties, we show below that the partition functions
are \emph{not equal} in the two formalisms. This result follows from an explicit
computation of the partition functions for two free bosons in the neighborhood of a two tori
degeneration points where $\Omega_{12}=0$. It then follows that there is
likewise no equality between the partition functions in the $\epsilon$- and 
$\rho$-formalisms for a lattice VOA. 

\medskip 
As shown in Theorem~12 of \cite{MT2}, we may relate the $\epsilon$- and 
$\rho$-formalisms in certain open neighborhoods of the two tori degeneration point, where  
$\Omega_{12}=0$. 
In the $\epsilon$-formalism, the genus two Riemann surface is parameterized by the domain 
\begin{equation}
\mathcal{D}^{\epsilon }=\{(\tau _{1},\tau _{2},\epsilon )\in  \mathbb{H}_{1}%
\mathbb{\times H}_{1}\mathbb{\times C}\ |\ |\epsilon |<\frac{1}{4}%
D(q_{1})D(q_{2})\},  \label{Deps}
\end{equation}
with $q_{a}=\exp(2\pi i \tau_{a})$ and $D(q)$ as in \eqref{Dom1om2}. 
In this case the two torus degeneration is, by definition, given by $\epsilon\rightarrow 0$.
 In the $\rho $-formalism, the two torus degeneration is described by the limit  \eqref{Omrho_degen}. In order to understand this more precisely 
we introduce the domain (\cite{MT2})
\begin{equation}
\mathcal{D}^{\chi }=\{(\tau ,w,\chi )\in \mathbb{H}_{1}\times \mathbb{%
C\times C}\ |\ (\tau ,w,-w^{2}\chi )\in \mathcal{D}^{\rho },0<|\chi |<\frac{1%
}{4}\},  \label{Dchidomain}
\end{equation}%
for $\mathcal{D}^{\rho }$ of \eqref{Drho} and $\chi=-\frac{\rho}{w^2}$ of \eqref{eq:chi2}.
The period matrix is determined by a $\Gamma _{1}$-equivariant holomorphic map 
\begin{eqnarray}
F^{\chi }:\mathcal{D}^{\chi } &\rightarrow &\mathbb{H}_{2},  \notag \\
(\tau ,w,\chi ) &\mapsto &\Omega ^{(2)}(\tau ,w,-w^{2}\chi ).  \label{Fchi}
\end{eqnarray}%
Then 
\begin{equation}
\mathcal{D}_{0}^{\chi }=\{(\tau ,0,\chi )\in \mathbb{H}_{1}\times \mathbb{C\times C}\ |\ 0<|\chi |<\frac{1}{4}\},
\label{D0chidomain}
\end{equation}%
is the space of two-tori degeneration limit points of the domain $\mathcal{D}^{\chi }$. 
We may  compare the two parameterizations on certain $\Gamma
_{1}$-invariant neighborhoods of a two tori degeneration point in both
parameterizations to obtain:
\begin{theorem}
\label{Theorem_epsrho_11map} (op. cit., Theorem 12) There exists a 1-1
holomorphic mapping between $\Gamma _{1}$-invariant open domains $\mathcal{I}%
^{\chi }\subset (\mathcal{D}^{\chi }\cup \mathcal{D}_{0}^{\chi })$ and $%
\mathcal{I}^{\epsilon }\subset \mathcal{D}^{\epsilon }$ where $\mathcal{I}%
^{\chi }$ and $\mathcal{I}^{\epsilon }$ are open neighborhoods of a two tori
degeneration point. $\hfill \square $
\end{theorem}

We next describe the explicit relationship between 
$(\tau _{1},\tau_{2},\epsilon )$ and $(\tau ,w,\chi )$ in more detail. Firstly, from Theorem~4 of \cite{MT2} we obtain%
\begin{eqnarray*}
\tpi\Omega _{11} &=&\tpi \tau _{{1}}+E_{2}\left( \tau _{2}\right) {%
\epsilon }^{2}+E_{2}\left( \tau _{1}\right) E_{2}\left( \tau _{2}\right) ^{2}%
{\epsilon }^{4}+O(\epsilon ^{6}), \\
\tpi\Omega _{12} &=&-{\epsilon }-E_{2}\left( \tau _{1}\right)
E_{2}\left( \tau _{2}\right) {\epsilon }^{3}+O(\epsilon ^{5}), \\
\tpi\Omega _{22} &=&\tpi \tau _{{2}}+E_{2}\left( \tau _{1}\right) {%
\epsilon }^{2}+E_{2}\left( \tau _{1}\right) ^{2}E_{2}\left( \tau _{2}\right) 
{\epsilon }^{4}+O(\epsilon ^{6}).
\end{eqnarray*}%
Making use of the identity
\begin{equation}
\frac{1}{\tpi }\frac{d }{d \tau }E_{2}(\tau )=5E_{4}(\tau )-E_{2}(\tau
)^{2},  \label{qdqE2}
\end{equation}%
it is straightforward to invert $\Omega _{ij}(\tau _{1},\tau _{2},\epsilon )$
to find
\begin{lemma}
\label{Lemma_Omeps}
In the neighborhood of the two tori degeneration point 
$r=\tpi\Omega _{12}=0$ of $\Omega \in \mathbb{H}_{2}$ we have 
\begin{eqnarray*}
\tpi \tau _{{1}} &=&\tpi\Omega _{11}
-E_{2}(\Omega_{22}) r^{2}
+5E_{2}(\Omega _{11}) E_{4}(\Omega_{22}) r^{4}+O(r^{6}), \\
{\epsilon } &=&-r+E_{2}(\Omega _{11})E_{2}(\Omega_{22})r^3+O(r^{5}), \\
\tpi \tau _{2} &=&\tpi\Omega_{22}
-E_{2}(\Omega _{11}) r^{2}
+5E_{2}(\Omega _{22}) E_{4}(\Omega _{11})r^4+O(r^6).
\end{eqnarray*}
$\hfill \Box$
\end{lemma}

From Theorem~\ref{Theorem_period_rho} we may also determine 
$\Omega_{ij}(\tau ,w,\chi )$ to $O(w^{4})$ in a
neighborhood of a two tori degeneration point to 
find
\begin{proposition}
\label{PropOmrhodegen}
For $(\tau ,w,\chi )\in \mathcal{D}^{\chi }\cup 
\mathcal{D}_{0}^{\chi }$ we have 
\begin{eqnarray*}
2\pi i\Omega _{11} &=&2\pi i\tau +
G(\chi )\sigma^{2}+G(\chi )^{2}E_{2}( \tau ) \sigma^{4}+O(w^{6}),  
\\
2\pi i\Omega _{22} &=&\log f(\chi )+E_{2}(\tau )\sigma^{2} +
\left( G(\chi )E_{2}( \tau )^{2}+\frac{1}{2}E_{4}( \tau ) \right) \sigma^{4}+O(w^{6})
\\
2\pi i\Omega _{12} &=&\sigma+
G(\chi )E_{2}(\tau)\sigma^{3}+O(w^{5}),
\end{eqnarray*}%
where $\sigma=w\sqrt{1-4\chi}$, $G(\chi )=\frac{1}{12}+E_{2}(q=f(\chi ))=O(\chi)$ and $f(\chi )$ is the Catalan series (\ref{Catalan}).
\end{proposition}
This result is an extension of   \cite{MT2}, Proposition~13 and the general proof proceeds along the same lines. 
For our purposes, it is sufficient to expand the non-logarithmic terms to $O(w^4,\chi^{0})$. 
Since $R(k,l)=O(\chi)$ and $d_{a}(k)=O(\chi^{1/2})$ then Theorem~\ref{Theorem_period_rho}  implies
\begin{eqnarray}
2\pi i\Omega _{11} &=&
2\pi i\tau +O(\chi),   \label{Om11_deg}\\
2\pi i\Omega _{22} &=& \log \chi +2\sum_{k\ge 2}\frac{1}{k}E_{k}(\tau)w^k +O(\chi), \label{Om22_deg}\\
2\pi i\Omega _{12} &=& w +O(\chi), \label{Om12_deg}
\end{eqnarray}
to all orders in $w$. In particular, we can readily
confirm Proposition~\ref{PropOmrhodegen} to $O(w^4,\chi^0)$.
Substituting \eqref{Om11_deg}--\eqref{Om12_deg} into Lemma~\ref{Lemma_Omeps}
and using \eqref{qdqE2} and \eqref{A.5} we obtain
\begin{proposition}
\label{Proposition_eps_rho_w4} For $(\tau ,w,\chi )\in \mathcal{D}^{\chi
}\cup \mathcal{D}_{0}^{\chi }$ we have
\begin{eqnarray*}
2\pi i\tau _{1} &=&\tpi \tau +\frac{1}{12}w^{2}+
\frac{1}{144} E_{2}( \tau )w^{4}+O(w^{6},\chi), \\
\tpi \tau _{2} &=&\log (\chi )
+\frac{1}{12}E_{4}( \tau ) w^{4}+O(w^{6},\chi), \\
\epsilon &=&-w-\frac{1}{12}E_{2}(\tau )w^{3}+O(w^{5},\chi).
\end{eqnarray*}
$\hfill \Box$
\end{proposition}

Define the ratio 
\begin{equation}
T_{\epsilon,\rho}(\tau ,w,\chi)=\frac{Z_{M^{2},\epsilon }^{(2)}(\tau _{1},\tau _{2},\epsilon )}
{Z_{M^{2},\rho }^{(2)}(\tau ,w,-w^2\chi )},
\label{Frat}
\end{equation}
for $\tau _{1},\tau _{2},\epsilon$ as given in Proposition~\ref{Proposition_eps_rho_w4}. 
From Theorems~8 of \cite{MT4} and Theorems \ref{Theorem_Z2_rho_G1} and \ref{Theorem_epsrho_11map} above we see that $T_{\epsilon,\rho}$ is $\Gamma_{1}$-invariant.  
From
Theorem \ref{Theorem_Degen_pt} for  $V=M^2$, we find in the two tori degeneration limit that
\begin{eqnarray*}
\lim_{w\rightarrow 0}T_{\epsilon,\rho}(\tau ,w,\chi)
&=&
f(\chi )^{-1/12},
\end{eqnarray*}
i.e., the two partition functions do not even agree in this limit!
The origin of this discrepancy may be thought to arise from the central charge dependent factors of $q^{-c/24}$ and $q_1^{-c/24}q_2^{-c/24}$ present in the definitions of $Z_{V,\rho }^{(2)}$ and $Z_{V,\epsilon }^{(2)}$ respectively (which, of course, are necessary for any modular invariance). 
One modification of the definition of the genus two partition functions compatible with the two tori degeneration limit might be: 
\begin{equation*}
Z_{V,\epsilon}^{\mathrm{new}(2)}(\tau _{1},\tau _{2},\epsilon )
=
\epsilon ^{-c/12}Z_{V,\epsilon}^{(2)}(\tau _{1},\tau _{2},\epsilon ),
\quad
Z_{V,\rho }^{\mathrm{new}(2)}(\tau ,w,\rho )=
\rho^{-c/24}Z_{V,\rho }^{(2)}(\tau ,w,\rho ).
\end{equation*} 
However, for $V=M^{2}$, we immediately observe that the ratio cannot be unity due to the incompatible $\Gamma_{1}$ actions arising from
\begin{eqnarray*}
\epsilon\rightarrow \frac{\epsilon}{c_{1}\tau_{1}+d_{1}},\quad
\rho\rightarrow \frac{\rho}{(c_{1}\tau+d_{1})^2}, 
\end{eqnarray*} 
as given in Lemmas~8 and 15 of \cite{MT2} (cf.\ \eqref{gam1_rho}).

\medskip
Consider instead a further $\Gamma _{1}$-invariant factor of $f(\chi )^{-c/24}$ in the definition of the genus two partition function in the 
$\rho $-formalism. Once again, we find that the partition functions do not agree in the neighborhood of a
two-tori degeneration point: 
\begin{proposition}
\label{Proposition_Zeps_rho_mod}
\begin{equation*}
f(\chi )^{1/12}T_{\epsilon,\rho}(\tau ,w,\chi)=
%\frac{Z_{M^{2},\epsilon }^{(2)}(\tau _{1},\tau _{2},\epsilon )}{f(\chi )^{-1/12}Z_{M^{2},\rho }^{(2)}(\tau ,w,-w^2\chi )}
1
-\frac{1}{288}E_{4}(\tau )w^{4}+O(w^{6},\chi).
\end{equation*}
\end{proposition}

\noindent \textbf{Proof.} 
As noted earlier, $R(k,l)=O(\chi)$ so that we immediately obtain 
\begin{equation*}
f(\chi )^{-1/12}Z_{M^{2},\rho }^{(2)}(\tau ,w,-w^2\chi )=
\frac{1}{{\eta (\tau)^{2}\eta (f(\chi))^{2}}}+O(\chi),
\end{equation*}
to all orders in $w$.
On the other hand,
$Z_{M^{2},\epsilon }^{(2)}(\tau _{1},\tau
_{2},\epsilon )$ of Theorem~5 of \cite{MT4} to $O(\epsilon ^{4})$
is given by
\begin{equation}
\frac{1}{\eta (\tau _{1})^{2}\eta (\tau _{2})^{2}}
\left[
1+E_{2}(\tau _{1})E_{2}(\tau _{2})\epsilon ^{2}+\left( E_{2}(\tau
_{1})^{2}E_{2}(\tau_{2})^{2}+15E_{4}(\tau _{1})E_{4}(\tau _{2})\right)
\epsilon^{4}
\right]+O(\epsilon ^{6}).
\label{Zeps4}
\end{equation}%
We expand this to $O(w^{4},\chi)$ using Proposition \ref{Proposition_eps_rho_w4}, \eqref{qdqE2} and 
\begin{equation*}
\frac{1}{2\pi i}\frac{d}{d \tau }\eta(\tau )=-\frac{1}{2}
E_{2}(\tau )\eta(\tau ),
\end{equation*}%
to eventually find that
\begin{gather*}
Z_{M^{2},\epsilon }^{(2)}(\tau _{1},\tau _{2},\epsilon )=
\frac{1}{\eta (q)^{2}\eta (f(\chi ))^{2}}
\left(
1-\frac{1}{288}E_{4}(\tau)w^{4}
\right)+O(w^{6},\chi). 
\end{gather*}%
$\hfill \Box$

\section{Final Remarks}
\label{Sect_Remarks}
Let us briefly and heuristically sketch how our
results compare to some related ideas in the physics and mathematics
literature. There is a wealth of literature concerning the bosonic string
e.g. \cite{GSW}, \cite{P}. In particular, the conformal anomaly implies that
the physically defined path integral partition function $Z_{\mathrm{string}}$
cannot be reduced to an integral over the moduli space $\mathcal{M}_{g}$ of
a Riemann surface of genus $g$ except for the 26 dimensional critical string
where the anomaly vanishes. Furthermore, for the critical string, Belavin
and Knizhnik argue that%
\begin{equation*}
Z_{\mathrm{string}}=\int\limits_{\mathcal{M}_{g}}\left\vert F\right\vert
^{2}d\mu ,
\end{equation*}%
where $d\mu $ denotes a natural volume form on $\mathcal{M}_{g}$ and $F$ is
holomorphic and non-vanishing on $\mathcal{M}_{g}$ \cite{BK}, \cite{Kn}.
They also claim that for $g\geq 2$, $F$ is a global section for the line
bundle $K\otimes \lambda ^{-13}$ (where $K$ is the canonical bundle and 
$\lambda $ the Hodge bundle) on $\mathcal{M}_{g}$ which is trivial by
Mumford's theorem \cite{Mu2}. In this identification,
the $\lambda ^{-13}$ section is associated with 26 bosons, the $K$ section
with a $c=-26$ ghost system and the vanishing conformal anomaly to the
vanishing first Chern class for $K\otimes \lambda ^{-13}$ \cite{N}.
More recently, some of these ideas have also been rigorously proved for a zeta
function regularized determinant of an appropriate Laplacian operator 
$\Delta _{n}$ \cite{McIT}. 
The genus two partition functions $%
Z_{M^{2},\epsilon }^{(2)}(\tau _{1},\tau _{2},\epsilon )$\ and $%
Z_{M^{2},\rho }^{(2)}(\tau ,w,\rho )$ constructed in \cite{MT4} and the
present paper for a rank $%
2 $ Heisenberg VOA should correspond in these approaches to a local
description of the holomorphic part of $\left( \frac{\det^{\prime }\Delta
_{1}}{\det N_{1}}\right) ^{-1}$ of refs. \cite{Kn}, $\ $\cite{McIT}, giving a
local section of the line bundle $\lambda ^{-1}$. 
Given these assumptions, it follows that $T_{\epsilon ,\rho }=Z_{M^{2},\epsilon
}^{(2)}/Z_{M^{2},\rho }^{(2)}\neq 1$ in the neighborhood of a two-tori
degeneration point where the ratio of the two sections is a non-trivial
transition function $T_{\epsilon ,\rho }$.

\medskip
In the case of a general rational conformal field theory, the conformal
anomaly continues to obstruct the existence of a global partition function
on moduli space for $g\geq 2$.
However, \emph{all} CFTs of a given central charge $c$ are believed to share
the same conformal anomaly e.g. \cite{FS}. Thus, the identification of the
normalized lattice partition and $n$-point functions of Section~\ref{Sect_LVOA}
reflect the equality of the first Chern class of some bundle associated to
a rank $c$ lattice VOA to that for $\lambda^{-c}$ with transition function $T_{\epsilon ,\rho }^{c/2}$. It is interesting to note that even in the case
of a unimodular lattice VOA with a unique conformal block (\cite{MS}, \cite%
{TUY}) the genus two partition function can therefore only be described
locally. It would obviously be extremely valuable to find a rigorous
description of the relationship between the VOA approach described here and
these related ideas in conformal field theory and algebraic geometry.

\section{Appendix}
\label{Subsect_Corr}
We list here some corrections to \cite{MT1} and \cite{MT2} that we needed above. 
\begin{description}
\item[(a)] 
Display (27) of \cite{MT1} should read 
\begin{equation}
\epsilon (\alpha ,-\alpha )=\epsilon (\alpha ,\alpha )=(-1)^{(\alpha ,\alpha
)/2}.  \label{A.1}
\end{equation}

\item[(b)] 
Display (45) of \cite{MT1} should read 
\begin{equation}
\gamma (\Xi )=(a,\delta _{r,1}\beta +C(r,0,\tau )\alpha _{k}+
\sum_{l\neq k}D(r,0,z_{kl},\tau )\alpha_{l}).  \label{A.2}
\end{equation}

\item[(c)] 
As a result of (a), displays (79) and (80) of \cite{MT1} are modified and now read 
\begin{eqnarray}
F_{N}(e^{\alpha },z_{1};e^{-\alpha },z_{2};q) &=&\epsilon (\alpha ,-\alpha )%
\frac{q^{(\beta ,\beta )/2}}{\eta ^{l}(\tau )}\frac{\exp ((\beta ,\alpha
)z_{12})}{K(z_{12},\tau )^{(\alpha ,\alpha )}},  \label{A.3} \\
F_{V_{L}}(e^{\alpha },z_{1};e^{-\alpha },z_{2};q) &=&\epsilon (\alpha
,-\alpha )\frac{1}{\eta ^{l}(\tau )}\frac{\Theta _{\alpha ,L}(\tau
,z_{12}/\tpi )}{K(z_{12},\tau )^{(\alpha ,\alpha )}}.  \label{A.4}
\end{eqnarray}

\item[(d)] 
The expression for $\epsilon(\tau,w,\chi)$ of display (172) of \cite{MT2} should read
\begin{equation}
\epsilon (\tau ,w,\chi ) =-w\sqrt{1-4\chi }(
1+\frac{1}{24}w^{2}E_{2}(\tau )(1-4\chi)+O(w^{4}))
\label{A.5}
\end{equation}
\end{description}

\pagebreak

\end{document}